\documentclass[12pt,reqno]{amsart}
\usepackage{amssymb}
\usepackage{cases}
\usepackage{bbm}
\usepackage{mathrsfs}
\usepackage{graphicx}
\usepackage{amsfonts}
\usepackage{enumerate}

\usepackage{ulem}
\usepackage{amssymb, amsmath, amsfonts, latexsym}
\usepackage{enumerate}
\usepackage{color}

\newtheorem{thm}{Theorem}[section]
\newtheorem{cor}[thm]{Corollary}
\newtheorem{lem}[thm]{Lemma}
\newtheorem{prop}[thm]{Proposition}
\newtheorem{example}[thm]{Example}
\newtheorem{remarks}[thm]{Remark}
\newtheorem{defn}[thm]{Definition}
\newtheorem{hyp}[thm]{Hypothesis}
\numberwithin{equation}{section}

\date{}

\newcommand{\ee}{\mathbb{E}}

\newcommand{\nn}{\mathbb{N}}
\newcommand{\rr}{\mathbb{R}}
\newcommand{\pp}{\mathbb{P}}
\newcommand{\qq}{\mathbb{Q}}

\def\FF{\mathcal F}

\def\LL{\mathcal L}

\def\NN{\mathcal N}

\def\vep{\varepsilon}
\def\<{\langle}
\def\>{\rangle}

\def\d"{^{\prime\prime}}

\def\bhyp{\begin{hyp}}
\def\nhyp{\end{hyp}}

\def\bbeq{\begin{equation}}
\def\nneq{\end{equation}}

\def\bdef{\begin{defn}}
\def\ndef{\end{defn}}

\def\bthm{\begin{thm}}
\def\nthm{\end{thm}}

\def\bprop{\begin{prop}}
\def\nprop{\end{prop}}

\def\brmk{\begin{remarks}}
\def\nrmk{\end{remarks}}

\def\bexa{\begin{example}}
\def\nexa{\end{example}}

\def\blem{\begin{lem}}
\def\nlem{\end{lem}}

\def\bcor{\begin{cor}}
\def\ncor{\end{cor}}

\def\bexe{\begin{exe}}
\def\nexe{\end{exe}}

\def\bprf{\begin{proof}}
\def\nprf{\end{proof}}

\def\dsp{\displaystyle}

\def\bdes{\begin{description}}
\def\ndes{\end{description}}

\title[Long-time behaviors of mean field interacting particle systems]
{Long-time behaviors of mean-field interacting particle systems related to McKean-Vlasov equations}

\author{Wei Liu}
\address{Wei Liu, School of Mathematics and Statistics, Wuhan University, Wuhan, Hubei 430072, PR China;
Hubei Key Laboratory of Computational Science, Wuhan University, Wuhan, Hubei 430072, PR China.}
\thanks{The first author is supported by NSFC (11731009), the Fundamental Research Funds for the Central Universities
2042020kf0217 and 2042020kf0031, and CSC}
\email{wliu.math@whu.edu.cn}

\author{Liming Wu}
\address{Liming Wu. Laboratoire de Math\'ematiques Blaise Pascal, CNRS-UMR
6620, Universit\'e Clermont-Auvergne (UCA), Campus Universitaire des Cezeaux, 3 Place Vasarely, 63178 Aubi\`ere, France.}
\email{Li-Ming.Wu@math.univ-bpclermont.fr}

\author{Chaoen Zhang}
\address{Chaoen Zhang. Laboratoire de Math\'ematiques Blaise Pascal, CNRS-UMR
6620, Universit\'e Clermont-Auvergne (UCA), Campus Universitaire des Cezeaux, 3 Place Vasarely, 63178 Aubi\`ere, France.}
\email{chaoenz@gmail.com}

\date{\today}
\begin{document}

\begin{abstract} In this paper, we investigate gradient estimate of the Poisson equation and the exponential convergence in the Wasserstein metric $W_{1,d_{l^1}}$, uniform in the number of particles, and uniform-in-time propagation of chaos for the mean-field weakly interacting particle system related to McKean-Vlasov equation. By means of the known approximate componentwise reflection coupling and with the help of some new cost function, we obtain explicit estimates for those three problems, avoiding the technical conditions in the known results. Our results apply when the confinement potential $V$ has many wells, the interaction potential $W$ has bounded second mixed derivative $\nabla^2_{xy}W$ which should be not too big so that there is no phase transition. As an application, we obtain the concentration inequality of the mean-field interacting particle system with explicit and sharp constants, uniform in time. Several examples are provided to illustrate the results.
\end{abstract}
\maketitle


\section{Introduction}
In this paper, we consider the following nonlinear McKean-Vlasov equation with initial condition $u_0$
\bbeq\label{MVE}
 \partial_t u_t=\nabla \cdot [\nabla u_t + u_t\nabla V + u_t(\nabla_x W\circledast u_t)],
\nneq
where the unknown $u_t$ is a time dependent probability density on $\rr^d$ ($d\ge1$), $V:\rr^d\to \rr$ is a
confinement potential and $W:\rr^d\times \rr^d  \to \rr$ is an interaction potential.
Here $\nabla$ and $\nabla\ \cdot$ (applied to a vector field) denote the gradient operator and the divergence operator respectively, while
$\nabla_x W$ stands for the gradient of $W$ with respect to (w.r.t. in short) the first variable, and
$$ \nabla_x W\circledast u_t(x):=\int_{\rr^d} \nabla_x W(x,y)u_t(y)dy.$$
When $W(x,y)=W_0(x-y)$ for some even potential $W_0$, $\nabla_x W\circledast u=\nabla W_0\ast u$ (the usual convolution).

The probabilistic equivalent version of (\ref{MVE}) is the following self-interacting stochastic differential equation (SDE in short):
\bbeq\label{SDE}\begin{cases}
dX_t=\sqrt{2}dB_t-\nabla V(X_t)dt - \nabla_x W\circledast \mu_t(X_t)dt,\\
X_0 \overset{law}{=} u_0(x)dx,
\end{cases}
\nneq
where $\mu_t$ is the law of $X_t$. The density $u_t$ of the law $\mu_t$ of $X_t$ at time $t$
is the solution of the McKean-Vlasov equation (\ref{MVE}) and {\it vice versa}. The existence and uniqueness of the solution of the SDE (\ref{SDE}) and the McKean-Vlasov equation (\ref{MVE}) have been extensively studied. The reader is referred to \cite{[MK],[G],[MS],[SZ]} and recent works \cite{[MV],[BLPR],[HSS]} as well as the references therein. For the convergence to equilibrium of solution $\mu_t$ as $t\to+\infty$, it is worth mentioning that Carrillo, McCann and Villani \cite{[CMV1]} obtained the explicit exponential convergence in entropy under various kinds of convexity conditions on the potentials $V$ and $W$, via their enlightening idea of interpreting the McKean-Vlasov equation as the gradient descent flow of the free energy on the space of probability measures equipped with the $L^2$-Wasserstein metric. Eberle  et al. \cite{[EGZ]} got the quantitative bounds on the exponential convergence in some appropriate transport cost to equilibrium for McKean-Vlasov equations by using Lyapunov condition and reflection coupling. Eberle \cite{[EA]} showed the exponential contractivity for diffusion semigroups w.r.t. Kantorovich distance by using componentwise reflection coupling methods and choosing appropriate distance functions. The reader is referred also to Luo and Wang \cite{[LW]} for the exponential convergence of diffusion semigroups w.r.t. the $L^p$-Wasserstein distance for all $p\ge 1$.

The McKean-Vlasov equation (\ref{MVE}) or (\ref{SDE}) is the idealization of the following interacting particle system of mean-field type when the number $N$ of particles goes to infinity:
\begin{equation}\label{MFS}\begin{cases}
dX_t^{i,N} = \sqrt{2} dB_t^i-\nabla V(X_t^{i,N})dt - \frac 1{N-1} \sum\limits_{j: j\ne i, 1\le j \le N}\nabla_x W(X^{i,N}_t, X^{j,N}_t) dt,\\
X_0^{i,N}=X_0^i, \ i=1,\cdots, N,
\end{cases}
\end{equation}
where the initial values $X_0^1, \cdots, X_0^N$ are i.i.d. random variables with common law $\mu_0(dx)=u_0(x)dx$, and $B_t^1 \cdots, B_t^N$
 are $N$ independent Brownian motions taking values in $\mathbb{R}^{d}$, independent of $X^{i}_0, 1\le i\le N$. In fact this is the goal of the so-called propagation of chaos: when the number $N$ of particles goes to infinity, the empirical measures $\frac 1N \sum_{i=1}^N \delta_{X^{i,N}_t}$ of the particle system (\ref{MFS}) (or the law of a single particle) converge weakly to the solution $\mu_t$ of the self-interacting diffusion (\ref{SDE}).

The propagation of chaos for the mean-field interacting particle systems has been widely studied during the last forty years. The early studies were concentrated on the propagation of chaos in bounded time intervals, see \cite{[Kac],[SZ],[MS]} and the references therein. The study on the propagation of chaos in the whole time interval $\rr^+$ is much more difficult and recent.
When the confinement potential $V$ is strictly convex and the interaction potential $W(x,y)=W_0(x-y)$ with $W_0$ strictly convex, Malrieu \cite{[MF1]} showed the uniform in time propagation of chaos by applying the logarithmic Sobolev inequality.  In the case that there is no confinement (i.e. $V\equiv 0$) and the interaction potential $W_0$ is strictly convex, Benachour et al. \cite{[BRTV1],[BRTV2]} proved propagation of chaos (but not uniform in time) and polynomial convergence to equilibrium; Malrieu \cite{[MF]} obtained the uniform in time propagation of chaos and exponential convergence to equilibrium for the particle system viewed from the center, by using  functional inequalities. When $W_0$ is degenerately convex and $V=0$, Cattiaux et al. \cite{[CGM]} showed the uniform in time propagation of chaos and exponential convergence to equilibrium by using synchronous coupling.

Without the convexity of $V$ and $W_0$, recently Durmus, Eberle, Guillin and Zimmer \cite{[DEGZ]} use the componentwise reflection coupling introduced in \cite{[EA]} to prove the exponential convergence in some Wasserstein metric and uniform in time propagation of chaos for weakly interacting mean-field particle system. For more results about propagation of chaos, we refer the reader to \cite{[DT],[HM],[JW],[L],[MM],[MMW]} and the references therein.

The main purpose of this paper is to investigate the exponential convergence in $L^1$-Wasserstein metric in the purpose of refining the previous results in \cite{[EA], [DEGZ]}, the concentration inequalities and the propagation of chaos of the mean-field weakly interacting particle system. Although we use the same approximate componentwise reflection coupling (\cite{[EA]}), our next approach will be quite different from  \cite{[EA],[DEGZ]}:

\begin{enumerate}
\item our starting point is some explicit gradient estimate of the Poisson equation, which implies moreover the concentration inequalities of the empirical mean of the interacting particle system, useful for numerical computation of solution $\mu_t$ of the McKean-Vlasov equation;

\item we will choose a different metric from that in \cite{[EA],[DEGZ]}, which allows us to obtain some explicit and almost sharp estimate of the exponential rate in the convergence of the interacting particles system to its equilibrium in the $W_1-$metric, uniform in the number $N$ of the particles.

\item
 As a by-product, we obtain some explicit estimate on the propagation of chaos, uniform in time.
\end{enumerate}

The paper is organized as follows. In the next section, we will present our framework and main results. The proofs are provided in Section 3 and section 4. The applications to concentration inequalities  are given in the last section.

\section{Main resutls}

\subsection{Framework: notations and conditions}

\subsubsection{Conditions on the dissipativity rate of a single particle}
First we introduce the dissipative rate $b_0(r)$ of the drift of one single particle in (\ref{MFS}) at distance $r>0$,
\bbeq \label{b0}\langle x-y, -[\nabla V(x)-\nabla V(y)]-[\nabla_x W(x,z)- \nabla_x W(y,z)]\rangle \le b_0(r)|x-y|
\nneq
holds for any $x,y,z\in\rr^d$ with $|x-y|=r$.
Throughout this paper we assume that $b_0(r)$ is a continuous function on $(0, +\infty)$ satisfying
\bbeq\label{b1} \limsup_{r\to +\infty} \frac{b_0(r)}{r} < 0,
\nneq
i.e. the drift of one particle is dissipative at infinity.

We also assume that
\bbeq\label{b2} \lim_{r\to0+}b^+_0(r) =0.
\nneq

Next we introduce an important reference function $h$ which enables us to obtain some new results, avoiding the technical parameters in \cite{[EA],[DEGZ]}.
For any function $f\in C^2(0,+\infty)$ and $r>0$, let $\LL_{ref}$ be the generator defined by
\begin{equation}\label{L1}
\LL_{ref}f(r):=4f^{\prime\prime}(r) + b_0(r)f^{\prime}(r).
\end{equation}
Let $h:\rr^+\to \rr^+$ be the function determined by: $h(0)=0$ and
\begin{equation}\label{PE2}
h^{\prime}(r)=\frac14 \exp\left(-\frac14 \int_0^r b_0(s)ds\right) \int_r^{+\infty} s\cdot\exp\left(\frac14 \int_0^s b_0(u)du\right) ds.
\end{equation}
$h$ is a well defined $C^2$ function by the dissipative condition (\ref{b1}), and it is a solution of the one-dimensional Poisson equation
\begin{equation}\label{PE1}
\LL_{ref} h(r)=4h^{\prime\prime}(r)+b_0(r)h^{\prime}(r)=-r,\ r>0
\end{equation}
with $h(0)=0$. This function was used by the second named author \cite{[Wu09JFA]} for functional and isoperimetric inequalities on Riemmanian manifolds.

\subsubsection{Kantorovich-Wasserstein $W_1$-metric}

For the configuration space $(\rr^d)^N$, instead of the usual Euclidean metric, we will use the $l^1$-metric (generalized Hamming metric)
$$d_{l^1}(x,y)=\sum_{i=1}^N |x^i-y^i|,\ \ x=(x^1,\cdots,x^N),\ y=(y^1,\cdots,y^N)\in(\rr^d)^N.$$
We consider the Kantorovich-Wasserstein distance w.r.t. $d_{l^1}$ metric on $(\rr^d)^N$, i.e., for any two probability measures $\mu$ and $\nu$ on $(\rr^{d})^N$,
$$W_{1,d_{l^1}}(\mu, \nu)=\inf_{P\in\Pi(\mu,\nu)}\iint_{(\rr^{d})^N\times (\rr^{d})^N}d_{l^1}(x,y)P(dx,dy)$$
where $\Pi(\mu,\nu)$ is the set of all couplings of $\mu,\nu$, i.e. the set of all probability measures on $(\rr^{d})^N\times (\rr^{d})^N$ whose marginal distributions of $x$ and $y$ are respectively $\mu$ and $\nu$.

Notice that for a $C^1$-function $g$ on $(\rr^d)^N$, its Lipschitzian norm $\|g\|_{Lip(d_{l^1})}$ w.r.t. $d_{l^1}$ coincides with $\max_{1\le i\le N}\|\nabla_i g\|_\infty$ where $\nabla_i$ is the gradient w.r.t. $x_i$.
By Kantorovich-Rubinstein duality relation,
$$W_{d_{l^1}}(\mu, \nu)=\sup_{g\in C_b^1((\rr^d)^N) : \max_{1\le i\le N}\|\nabla_i g\|_\infty \le 1}\left(\int g d\mu -\int g d\nu\right)
$$
When $N=1$, we write simply $W_1$ for $W_{1,d_{l^1}}$.

We notice that {\it for two probability measures $\mu,\nu$ on $(\rr^d)^N$,
\bbeq\label{W1P}
\sum_{i=1}^N W_1(\mu^i,\nu^i) \le W_{d_{l^1}}(\mu, \nu)
\nneq
and the equality holds when $\dsp \mu=\otimes_{i=1}^N \mu^i, \nu=\otimes_{i=1}^N \nu^i$ are product measures, }
where $\mu^i$ (resp. $\nu^i$) is the marginal distribution of $x_i$ of $\mu$ (resp. $\nu$). In fact if $X=(X^1,\cdots, X^N), Y=(Y^1,\cdots, Y^N)$ are two random vectors such that the law of $(X,Y)$ is an optimal coupling of $(\mu,\nu)$ in $W_{1,d_{l^1}}$, then for each $i$, the law of $(X^i,Y^i)$ is a coupling of $(\mu^i,\nu^i)$, so
$$
W_{d_{l^1}}(\mu, \nu)= \ee d_{l^1}(X,Y) =\sum_{i=1}^N \ee |X^i-Y^i| \ge \sum_{i=1}^N W_1(\mu^i,\nu^i).
$$
When $\mu$, $\nu$ are product measures, let $(X^i, Y^i)$ (or its joint law) be an optimal coupling of $(\mu^i,\nu^i)$ for $W_1(\mu^i,\nu^i)$ so that $(X^1, Y^1),\cdots, (X^N, Y^N)$ are independent. Then $(X=(X^i)_{1\le i\le N}, Y=(Y^i)_{1\le i\le N}$) is a coupling of $(\mu,\nu)$, so we get
$$
\sum_{i=1}^N W_1(\mu^i,\nu^i)= \sum_{i=1}^N \ee |X^i-Y^i|= \ee d_{l^1}(X,Y)\ge W_{d_{l^1}}(\mu, \nu)
$$
i.e. the equality in (\ref{W1P}) holds in the prodcut measures case. (This is well known.)

\subsection{An explicit gradient estimate of the Poisson equation and its applications in concentration inequalities}
Let $\{P_t^{(N)}\}_{t\ge 0}$ be the transition semigroup of the mean-field interacting particle system (\ref{MFS}), whose generator is given by
$$
\LL^{(N)} f(x^1,\cdots,x^N) = \sum_{i=1}^N \left(\Delta_i f - \nabla V(x^i)\cdot\nabla_i f - \frac 1{N-1} \sum_{j\ne i} \nabla_x W(x^i,x^j)\cdot \nabla_i f\right).
$$
Its unique invariant probability measure is the mean-field Gibbs measure, given by
$$
\mu^{(N)}(dx^1,\cdots,dx^N) =\frac 1{C_N} \exp\left(-\sum_{i=1}^N V(x^i) - \frac 1{N-1} \sum_{1\le i<j\le N} W(x^i,x^j)\right)dx^1\cdots dx^N,
$$
where $C_N$ is the normalization constant.

We introduce the following key assumption on the interaction potential:
$$ \text{(\bf H)}:  \ \|\nabla^2_{xy}W\|_{\infty} \|h^{\prime}\|_{\infty} < 1$$
where $h$ is given by (\ref{PE2}), $\|h^{\prime}\|_{\infty}:=\sup_{r\ge 0}h^{\prime}(r)$, and $\nabla^2_{xy} W=(\frac{\partial^2}{\partial x_i\partial y_j}W)_{1\le i, j\le d}$,
$$\|\nabla^2_{xy}W\|_\infty:=\sup_{x,y\in\rr^d}\sup_{z\in\rr^d, |z|=1}|\nabla^2_{xy}W(x,y)z|.$$
Notice that when the dissipativity at infinity condition (\ref{b1}) is satisfied, $b_0(r)$ can be taken as $-c_1r+c_2$ (with $c_1,c_2>0$), so $\|h'\|_\infty:=\sup_{r\ge0}h'(r)<+\infty$.

This assumption is a translation of Dobrushin-Zegarlinski's uniqueness condition in the framework of mean field, and
it implies that the mean field has no phase of transition (see \cite{[GLWZ]}).

Notice that under the assumption $(\bf H)$ and (\ref{b1}), both the equations (\ref{SDE}) and (\ref{MFS}) have unique strong solutions. On the space of continuous paths $C([0,T], (\rr^d)^N)$ where $T\in (0,+\infty]$, we consider the $L^1$-metric
\bbeq\label{dL}
d_{L^1[0,T]}(\gamma_1, \gamma_2) :=\int_0^{T} d_{l^1}(\gamma_1(t), \gamma_2(t)) dt.
\nneq
Given the starting point $x\in (\rr^d)^N$, let $\pp_x$ be the law of $X^{(N)}=(X_t^{(N)})_{t\ge0}$ with $X_0^{(N)}=x$.

\begin{thm}\label{W1} Assume (\ref{b1}), (\ref{b2}) and $(\bf H)$. For any $x_0=(x_0^1,\cdots,x_0^N)\in(\rr^d)^N$
 and $y_0=(y_0^1,\cdots,y_0^N)\in(\rr^d)^N$,
 we have
\begin{equation}\label{W2}
\aligned
\int_0^{+\infty} W_{d_{l^1}}(P_t^{(N)}(x_0,\cdot), P_t^{(N)}(y_0,\cdot)) dt &\le W_{1,d_{L^1[0,\infty]}}(\pp_{x_0}, \pp_{y_0})\\
&\le \frac{1}{1-\|\nabla^2_{xy}W\|_{\infty} \|h^{\prime}\|_{\infty}} \sum_{i=1}^N h(|x_0^i-y_0^i|).
\endaligned\end{equation}
 In particular for any $g\in C_b^1((\rr^d)^N)$ with $\mu^{(N)}(g)=0$,  the solution $G$ of the Poisson equation $-\LL^{(N)} G=g$ with $\mu(G)=0$ satisfies
\begin{equation}\label{cor1}
\|\nabla_i G\|_{\infty} \le  c_{Lip} \cdot \max_{1\le j\le N}\|\nabla_j g\|_\infty,
1\le i \le N,
\end{equation}
where \bbeq\label{clip} c_{Lip}:=\frac{h^{\prime}(0)}{1-\|\nabla^2_{xy}W\|_{\infty} \|h^{\prime}\|_{\infty}}\nneq
and
$$h^{\prime}(0)=\frac14 \int_0^{+\infty} s\cdot\exp\left(\frac14 \int_0^s b_0(u)du\right)ds.$$
\end{thm}

By the theorem above we can immediately obtain the following result about the nonlinear McKean-Vlasov equation (\ref{MVE}).

\bcor\label{cor11}  Under the same assumptions as in Theorem \ref{W1},
for any two solutions $\mu_t,\nu_t$ of the self-interacting diffusion (\ref{SDE}) with the initial distributions $\mu_0,\nu_0$ with finite second moment respectively, we have
\begin{equation}\label{cor11a}
\int_0^\infty W_1(\mu_t, \nu_t) dt \le  \frac{\|h^\prime\|_\infty}{1-\|\nabla^2_{xy}W\|_\infty \|h^\prime\|_\infty} W_1(\mu_0,\nu_0).
\end{equation}
\ncor

\bprf By \eqref{W2} in Theorem \ref{W1} and the fact that $$h(r) \le h(0)+ \|h^\prime\|_\infty \cdot r=\|h^\prime\|_\infty \cdot r,\ \forall r\ge 0$$
 we have
\bbeq\label{cor11c}\int_0^\infty W_{1,d_{l^1}}(\mu_0^{\otimes N}P^{(N)}_t, \nu_0^{\otimes N} P^{(N)}_t)dt  \le \frac{\|h^\prime\|_\infty}{1-\|\nabla^2_{xy}W\|_\infty \|h^\prime\|_\infty}W_{1,d_{l^1}}(\mu_0^{\otimes N}, \nu_0^{\otimes N}).\nneq
Notice that $\mu^{(N)}_t:=\mu_0^{\otimes N}P^{(N)}_t$ and $\nu^{(N)}_t:=\nu_0^{\otimes N}P^{(N)}_t$ are symmetric probability measures on $(\rr^d)^N$ and their marginal distributions $\mu^{(i,N)}_t$, $\nu^{(i,N)}_t$ of $x_i$ converge weakly to $\mu_t, \nu_t$ (respectively) by the finite time propagation of chaos. By using (\ref{W1P}) we have
$$
N W_1(\mu^{(1,N)}_t, \nu^{(1,N)}_t) = \sum_{i=1}^N W_1(\mu^{(i,N)}_t, \nu^{(i,N)}_t)\le W_{1,d_{l^1}}(\mu^{(N)}_t, \nu^{(N)}_t)
$$
and then by the lower semi-continuity of $W_1$ in the weak convergence topology,
\bbeq\label{cor11b} W_1(\mu_t, \nu_t) \le \liminf_{N\to+\infty} W_1(\mu^{(1,N)}_t, \nu^{(1,N)}_t) \le \liminf_{N\to+\infty} \frac 1N W_{1,d_{l^1}}(\mu^{(N)}_t, \nu^{(N)}_t).
\nneq

Combining \eqref{cor11c} and \eqref{cor11b} together, we obtain by Fatou's lemma,
$$\aligned \int_0^\infty W_1(\mu_t, \nu_t) dt &\le \liminf_{N\to+\infty} \frac 1N \int_0^\infty  W_{1,d_{l^1}}(\mu^{(N)}_t, \nu^{(N)}_t) dt\\
&\le \frac{\|h^\prime\|_\infty}{1-\|\nabla^2_{xy}W\|_\infty \|h^\prime\|_\infty} \liminf_{N\to+\infty} \frac 1N W_{1,d_{l^1}}(\mu_0^{\otimes N}, \nu_0^{\otimes N})\\
&=\frac{\|h^\prime\|_\infty}{1-\|\nabla^2_{xy}W\|_\infty \|h^\prime\|_\infty} W_1(\mu_0,\nu_0)
\endaligned$$
where the last equality follows by (\ref{W1P}). That completes the proof.
\nprf

As an application of Theorem \ref{W1} to the concentration inequality, we have the following result about the Gaussian concentration of the $U$-statistics, which is a straightforward application of a general result in Proposition \ref{C}. The proofs are given in the last section \S 5.

For any $1\le m\le N$, let $f_m: (\rr^d)^m\to \rr $ be a measurable and symmetric function. The $U$-statistic of order $m$ with kernel $f_m$ is defined by \bbeq U_N(f_m)(x^1,\cdots,x^N) = \frac 1{|I_N^m|} \sum_{(i_1,\cdots,i_m)\in I_N^m}  f_m(x^{i_1},\cdots, x^{i_m}),\ \forall (x^1,\cdots,x^N)\in (\rr^d)^N,\nneq
where
\bbeq\label{INm}
I_N^m:=\{(i_1,\cdots,i_m)\in \nn^k|i_1, \cdots, i_m \text{ are different }, 1\le i_1,\cdots, i_m\le N\}\nneq
and $|I_N^m|$ denotes the number of elements in $I_N^m$ (equal to $N !/(N-m) !$).

Next we introduce the following {\it Gaussian integrability} assumption of the initial distribution $\mu_0$:
\bbeq\label{Ua} \int_{\rr^d} e^{\lambda_0 |x|^2}\mu_0(dx)<+\infty,\ \text{for some}\ \lambda_0>0
\nneq
which is equivalent to say that there is some Gaussian concentration constant $c_{G}(\mu_0)>0$ such that
\bbeq\label{Ub}
\int_{\rr} e^{f(x)-\mu_0(f)} d\mu_0(x) \le \exp\left(\frac {c_G(\mu_0)}2 \|f\|_{Lip}^2\right)
\nneq
for all Lipschitzian functions $f$ on $\rr^d$ (w.r.t. the usual Euclidean distance).

\brmk\label{remc}{\rm The equivalence between the {\it Gaussian integrability} (\ref{Ua}) and the {\it Gaussian concentration inequality} (\ref{Ub}) was established by H. Djellout, A. Guillin and the second named author \cite{[DGW]}, and  (\ref{Ub}) is the famous characterization of Bobkov-G\"otze \cite{[BoG]} of the transport-entropy inequality.
By the tensorization of the transport-entropy inequality for product measure, (\ref{Ub}) implies that for any $N\ge 1$,

\bbeq\label{Uc}
\int_{(\rr^d)^N} e^{g(x)-\mu_0^{\otimes N}(g)} d\mu_0^{\otimes N}(x) \le \exp\left(\frac N2 c_{G}(\mu_0) \|g\|^2_{Lip(d_{l^1})}\right)
\nneq
for all Lipschitzian functions $g$ on $(\rr^d)^N$.
}
\nrmk

\bcor\label{U} Assume the conditions in Theorem \ref{W1} and the Gaussian integrability \eqref{Ua} of the initial distribution $\mu_0$. Let $f_m\in C^2((\rr^d)^m,\rr)$ be measurable, symmetric, and $1$-Lipschitz w.r.t. the $d_{l^1}$-metric on $(\rr^d)^m$, i.e. $\max_i \|\nabla_i f\|_\infty\le 1$.
Then for any $\lambda, T>0$, we have
\begin{equation}\begin{aligned}\label{cor3}
&\mathbb{E}\exp\left(\frac{\lambda}{T} \left[\int_0^T U_N(f_m)(X_t^{1,N},\cdots,X_t^{N,N}) dt - \int_0^T \mathbb{E}f_m(X_t^{1,N},\cdots,X_t^{m,N}) dt\right]\right)\\
&\le \exp\left(\frac{m^2\lambda^2c_{Lip}^2}{2NT}\left(1+\frac{c_G(\mu_0)}{T}\right)\right),\\
\end{aligned}\end{equation}
where $c_{Lip}$ is the same as given in (\ref{clip}).
In particular we have for any $\delta>0$

\begin{equation}\label{cor4}\begin{aligned}
&\mathbb{P}\left\{\frac{1}{T} \int_0^T U_N(f_m)(X_t^{1,N},\cdots,X_t^{N,N}) dt - \frac{1}{T} \int_0^T \mathbb{E} f_m(X_t^{1,N},\cdots,X_t^{m,N}) dt >\delta \right\}\\
&\le \exp\left(-\frac{(1-\|\nabla^2_{xy}W\|_{\infty} \|h^{\prime}\|_{\infty})^2}{2m^2(h^{\prime}(0))^2(1+c_G(\mu_0)/T)} NT\delta^2\right).\\
\end{aligned}\end{equation}
\ncor
The concentration inequality (\ref{cor4}) is sharp when $V$ is quadratic and $W=0$, see Example \ref{Gaussian} for explicit expression of all involved constants in the Gaussian case.

\subsection{Exponential convergence of the particle system in the $W_{1,d_{l^1}}$-metric}
\begin{thm}\label{W1c} Assume (\ref{b1}) and $(\bf H)$. Suppose that there exists a constant $M\in\rr$
such that
\bbeq\label{b5} b_0(r)\le rM, \forall r>0\nneq
(this condition is stronger than (\ref{b2})), then for any $\vep>0$ such that
\bbeq\label{b6} K_{\vep}:=\frac{1-\|\nabla^2_{xy}W\|_{\infty}\|h^{'}\|_{\infty}- \vep (M+\|\nabla^2_{xy}W\|_{\infty})}{\|h^{\prime}\|_{\infty}+\vep} >0,
\nneq
we have for any $x_0,y_0\in (\rr^d)^N$
\begin{equation}\label{W1c1}
W_{d_{l^1}}(P_t^{(N)}(x_0,\cdot), P_t^{(N)}(y_0,\cdot)) \le A_{\vep} e^{-K_{\vep}t} d_{l^1}(x_0, y_0),\ \forall t\ge 0,
\end{equation}
where
\bbeq\label{cKA}\aligned
A_{\vep}=\sup_{r>0}\frac{r}{h(r)+\vep r} \cdot \sup_{r>0}\frac{h(r)+\vep r}{r}.
\endaligned \nneq
\end{thm}

\begin{remarks}\label{W1cr} {\rm  An easy estimate of $A_{\vep}$ is $A_\vep\le \frac{\sup_{r\ge0}h^{'}(r) +\vep}{\inf_{r\ge0}h^{'}(r) +\vep}$ (since $h(0)=0$).
Note that when $M+\|\nabla^2_{xy}W\|_{\infty}>0$, the exponential rate $K_{\vep}$
increases (then better and better) as $\vep$ decreases to $0$, but $A_{\vep}$ may explode once if $\inf_{r\ge0}h^{'}(r)=0$.

}
\end{remarks}

\brmk\label{BE}{\rm Notice that (\ref{b5}) is equivalent to say that
$$
\nabla^2 V(x) + \nabla^2_{xx} W(x,y) \ge - M I, \ x,y\in\rr^d.
$$
When $\kappa:=-M - \|\nabla^2_{xy} W\|_\infty>0$, we see that the Hessian of the Hamiltonian
$$H(x^1,\cdots,x^N)=\sum_{i=1}^N V(x^i) + \frac 1{N-1}\sum_{1\le i<j\le N} W(x^i,x^j)$$
is bounded from below by $\kappa I$ (this estimate of the lower bound of the Bakry-Emery curvature is sharp if $\nabla^2_{xy}W$ is constant and definitely nonnegative). Notice that when $M<0$, we can take $b_0(r)=Mr$, so $h'(r)=-1/M$. Then $\kappa>0$ if and only if
{\bf (H)} is satisfied. The advantage of our condition {\bf (H)} (w.r.t. the positive curvature condition) is: it does not depend on the curvature but on the dissipativity, it holds even if $V$ has many wells once if the interaction is weak enough.

 In the positive curvature $\kappa>0$ case we have by Bakry-Emery's curvature characterization
$$W_1(P_t^{(N)}(x,\cdot), P_t^{(N)}(y,\cdot))\le e^{-\kappa t} |x-y|$$
in the Euclidean metric on $(\rr^d)^N$. On the other hand as above $b_0(r)=Mr, h(r)=-r/M$, we see that $K_\vep\to -M-\|\nabla^2_{xy}W\|_\infty=\kappa$ as $\vep\to +\infty$, and $A_\vep\equiv1$, so (\ref{W1c1})
 yields
\bbeq\label{W1cb}
W_{1,d_{l^1}}(P_t^{(N)}(x,\cdot), P_t^{(N)}(y,\cdot))\le e^{-\kappa t} d_{l^1}(x,y),
\nneq
a new but not at all surprising observation.
}\nrmk

Theorem \ref{W1c} above will give us an explicit exponential convergence in $W_1$ of the nonlinear McKean-Vlasov equation (\ref{MVE}). For the exponential convergence in entropy of the nonlinear McKean-Vlasov equation (\ref{MVE}) under the condition {\bf (H)}, see Guillin {\it et al.} \cite{[GLWZ]}.

\bcor\label{cor5}  Under the same assumptions as in Theorem \ref{W1c}, for any  $\vep>0$ so that $K_\vep>0$ (i.e. (\ref{b6})),
we have for the solutions $\mu_t,\nu_t$ of the self-interacting diffusion (\ref{SDE}) with the initial distributions $\mu_0,\nu_0$ which have finite second moments respectively,
\begin{equation}\label{cor5a}
W_{1}(\mu_t, \nu_t) \le  A_{\vep} e^{-K_{\vep}t} W_1(\mu_0,\nu_0),\ \forall t\ge 0,
\end{equation}
where $K_{\vep}$ and $A_{\vep}$ are given by \eqref{b6} and (\ref{cKA}) respectively.
\ncor
\bprf The proof of this corollary is similar to that of Corollary \ref{cor11}, and we utilize the same notations as in the Corollary \ref{cor11}.  First by Theorem \ref{W1c}, we have for any $t\ge 0$
$$W_{1,d_{l^1}}(\mu_0^{\otimes N}P^{(N)}_t, \nu_0^{\otimes N} P^{(N)}_t)  \le A_{\vep} e^{-K_{\vep}t}W_{1,d_{l^1}}(\mu_0^{\otimes N}, \nu_0^{\otimes N}).$$

Combining the inequality above with \eqref{cor11b}, we obtain
$$\aligned
W_1(\mu_t, \nu_t)&\le A_{\vep} e^{-K_{\vep}t} \liminf_{N\to+\infty} \frac 1N  W_{1,d_{l^1}}(\mu_0^{\otimes N}, \nu_0^{\otimes N})\\
&= A_{\vep} e^{-K_{\vep}t} W_1(\mu_0,\nu_0)
\endaligned$$
the desired result.
\nprf

\subsection{Propagation of chaos in large time}

We have the following uniform in time propagation of chaos.

\begin{thm}\label{Wp} Assume (\ref{b1}), (\ref{b5}) and $(\bf H)$. Suppose that there exist some positive constants $c_1,c_2,c_3$ such that
 \bbeq\label{MEaa} \langle x, \nabla V(x)\rangle\ge c_1|x|^2-c_2, \forall x\in\rr^d
 \nneq
 and
 \bbeq\label{MEbb} \langle z, \nabla^2_{xx}W(x,y) z\rangle\ge -c_3 |z|^2, \forall x,y,z\in\rr^d.
 \nneq
 Assume
\bbeq\label{MEcc}c_1-c_3-\|\nabla^2_{xy}W\|_{\infty}>0.\nneq

Then for any $\vep > 0$ such that $K_\vep >0$, and $\tilde{\vep} \in (0, c_1-c_3-\|\nabla^2_{xy}W\|_{\infty}),$
the following estimates of propagation of chaos hold for the mean-field interacting particle system (\ref{MFS}) with any initial probability measure $\mu_0$ having finite second moment:

\begin{enumerate}[(a)]

\item {\bf (path-type propagation of chaos)} for any $T>0$, $1\le k\le N$, denote $\pp_\nu(\cdot)=\int_{(\rr^d)^N} \pp_x(\cdot) d\nu(x)$ the law of $(X^{(N)}_t)_{t\ge0}$ with the initial distribution $\nu$, $\pp_\nu^{[1,k],N}|_{[0,T]}$ the joint law of paths of the $k$ particles $((X^{i,N}_t)_{t\in [0,T]}, 1\le i\le k)$ in time interval $[0,T]$,
and $\qq_{\mu_0}$ the law of
the self-interacting diffusion $(X_t)_{t\ge0}$ with the initial distribution $\mu_0$. We have
\bbeq\label{MEdd}
\frac 1k W_{1, d_{L^1[0,T]}} (\pp_{\mu_0^{\otimes N}}^{[1,k],N}|_{[0,T]}, \qq_{\mu_0}^{\otimes k}|_{[0,T]}) \le \frac{T}{\sqrt{N-1}}\frac{\|\nabla^2_{xy}W\|_{\infty}\|h^\prime\|_\infty }{1-\|\nabla^2_{xy}W\|_{\infty}\|h^{\prime}\|_{\infty}}\cdot \max\{m_2(\mu_0), \hat c(\vep)\}
\nneq
where \bbeq\label{mc}\aligned &m_2(\mu_0)=\left(\int_{\rr^d} |x|^2d\mu_0(x)\right)^{\frac 12},\\
&\hat c(\vep)= \left(\frac{d+c_2+ \frac 1{4\tilde \vep}|\nabla_x W(0,0)|^2}{c_1-c_3-\|\nabla^2_{xy}W\|_{\infty}-\tilde\vep}\right)^{\frac 12}.
\endaligned\nneq

\item {\bf (Uniform in time propagation of chaos)} for all time $t>0$ and any $1\le k\le N$:
\begin{equation}\label{CP0}
 W_{1, d_{l^1}}( \mu_t^{[1,k],N}, \mu^{\otimes k}_t) \le \frac{k}{\sqrt{N-1}}\frac{A_\vep}{K_\vep}\|\nabla^2_{xy}W\|_{\infty} \max\{m_2(\mu_0), \hat c(\vep)\}
\end{equation}where  $\mu_t=u_tdx$ is the solution of the McKean-Vlasov equation (\ref{MVE}), and $\mu_t^{[1,k],N}$ is the joint law of the $k$ particles $(X_t^{i,N}, 1\le i\le k)$ in the mean-field system (\ref{MFS}) of interacting particles $(X_t^{i,N})_{1\le i\le N}$ with  $X^{i,N}_0,
1\le i\le N$ i.i.d. of law $\mu_0$ (independent of $(B_t^{i,N})_{1\le i\le N, t\ge0}$), and the constants $K_\vep$, $A_\vep$, $m_2(\mu_0)$ and $\hat c(\vep)$  are given in \eqref{b6}, (\ref{cKA}) and \eqref{mc} respectively.

\end{enumerate}

\end{thm}

\brmk\label{remb} {\rm The time-uniform propagation of chaos is much more difficult than the bounded time propagation of chaos, accomplished in the 80-90's of the last century.  The physical reason is that the time-uniform propagation of chaos fails in the regime of phase transition. That is why we impose the condition {\bf (H)}, which excludes the phase transition.

The reader is referred to \cite{[BGM],[CGM],[DT],[DEGZ],[L]} and the references therein for recent studies and progresses on this subject. The main new point here is that our estimate (\ref{CP0}) is explicit and relatively neat.

}
\nrmk

\brmk\label{rema} {\rm All the results presented in this paper can be extended to more general case:
$$dX_t=\sqrt{2}dB_t+b(X_t,\mu_t)dt$$
where $\mu_t$ is the law of $X_t$, if $b$ satisfies some dissipative condition in $x$ (uniformly in $\mu$) and a Lipchitzian condition in $\mu$ with sufficiently small Lipschitzian constant. For the sake of clarity, we deal only with the case of $b(X_t,\mu_t)=-\nabla V(X_t)-\nabla_x W\circledast \mu_t(X_t)$ in this paper.
}
\nrmk

\subsection{Examples}

We first present the Gaussian model for which the constants in Theorem \ref{W1} and Theorem \ref{W1c} become exact, showing their sharpness.

\bexa\label{Gaussian} {\bf (Gaussian model)} Let $d=1$, and $$V(x)=\beta \frac{ x^2}{2},\ W(x,y)=-\beta K xy$$ where $\beta>0$ is the inverse temperature, $K\ge0$.

{\rm For this model, by some simple calculations we have
$$b_0(r)=-\beta r,\ \forall r>0.$$
and
$$h^{\prime}(r)\equiv\beta^{-1},\ \forall r\ge 0.$$
It is obvious that conditions (\ref{b1}) and (\ref{b2}) hold, and the assumption $(\bf H)$ holds once if
\bbeq\label{cw5} K < 1.\nneq
But this condition is equivalent to say that the matrix $A=(a_{ij})_{1\le i,j\le N}$ is positively definite, where
$$
a_{ii}= \beta, a_{ij}=\frac {-\beta K}{N-1}, i\ne j.
$$
$A$ must be the inverse of the covariance matrix of the Gaussian measure $\mu^{(N)}$. In other words {\bf (H)} is equivalent to well defining the equilibrium probability measure $\mu^{(N)}$.

Note that $\|\nabla^2_{xy}W\|_{\infty}=\beta K$, so we have $c_{Lip}=\frac{1}{\beta(1-K)}$  under (\ref{cw5}). Moreover (\ref{b5}) is satisfied with $M=-\beta$.

\medskip
$\bullet$ {\it Sharpness of Theorem \ref{W1}.}  The gradient estimate (\ref{cor1}) in
Theorem \ref{W1} tells us: if $-\LL^{(N)}G=g$, then
$$
\|\nabla_i G\|_\infty \le \frac{1}{\beta(1-K)} \max_{i}\|\nabla_i g\|_\infty.
$$
Let us show that it becomes equality for $g(x^1,\cdots, x^N)=\sum_{i=1}^N x^i$. In fact
$$
\LL^{(N)} g(x^1,\cdots, x^N) =-\sum_{i} \beta x^i + \sum_{i} \frac 1{N-1} \sum_{j\ne i} \beta K x^j =- \beta (1-K) g.
$$
In other words $G=\frac 1{\beta(1-K)} g$ for which the gradient estimate above becomes equality. As the gradient estimate (\ref{cor1}) comes from  (\ref{W2}), the process level $W_{1,d_{L^1}}$ estimate  (\ref{W2}) is sharp too.

 \medskip
$\bullet$ {\it Sharpness of Theorem \ref{W1c}.} As $\vep\to+\infty$ in (\ref{b6}), we have by Theorem \ref{W1c}
$$
W_{1,d_{l^1}}(P_t^{(N)}(x_0,\cdot), P_t^{(N)}(y_0,\cdot)) \le e^{-\beta(1-K)t} d_{l^1}(x_0,y_0).
$$
This is equivalent to say that
$$
\max\|\nabla_i P_t^{(N)} g\|_\infty \le e^{-\beta(1-K)t} \max\|\nabla_i g\|_\infty.
$$
But it becomes equality for $g=\sum_{i=1}^N x^i$ : in fact as $\LL^{(N)}g=-\beta(1-K) g$,
$$
P_t^{(N)} g = e^{-\beta(1-K)t} g.
$$
Hence the exponential convergence result \eqref{W1c1} in Theorem \ref{W1c} is sharp.

Of course for this Gaussian model all results in Theorems \ref{W1} and \ref{W1c} can be derived easily by using the synchronous coupling, or from the commutativity relation
$$
\nabla P_t^{(N)} g = e^{-At} P_t^{(N)}\nabla g
$$
which is one of the origins of the Bakry-Emery curvature.

}\nexa

Next we give another two typical models to illustrate our results.
\bexa\label{Curie-Weiss} {\bf (Curie-Weiss mean-field lattice model)} Let $d=1$, and $$V(x)=\beta(x^4/4 -x^2/2),\ W(x,y)=-\beta K xy$$ where $\beta=\frac 1{\kappa T}>0$ ($\kappa$ is the Boltzmann constant) is the inverse temperature, $K\in \rr^*$. This model is ferromagnetic or anti-ferromagnetic according to $K>0$ or $K<0$.

{\rm
By an elementary calculation, we get $$b_0(r)=\beta r(1-r^2/4),\ \forall r>0.$$
It is obvious that conditions (\ref{b1}) and (\ref{b2}) are satisfied and (\ref{b5}) holds with $M=\beta$.

For the assumption $(\bf H)$, first notice that $\|\nabla_{xy}^2 W\|_{\infty}=|K|\beta$. Next we estimate $\|h^\prime\|_{\infty}$. By (\ref{PE2}) and some calculations, we have for any $r\ge0$
$$\aligned
h^\prime(r)&=\frac14 \exp(\beta(r^4-8r^2)/64)\int_r^{+\infty}s\cdot \exp(\beta(8s^2-s^4)/64)ds\\
&=\frac14 e^{\beta/4} \exp(\beta(r^4-8r^2)/64) \int_{r^2/2}^{+\infty}\exp(-\beta(u-2)^2/16)du.\\
\endaligned$$
When $\frac{r^2}{2}>2$, i.e. $r>2$, we have
\bbeq\label{cw1}h^\prime(r)\le \frac14 e^{\beta/4} \exp(\beta(r^4-8r^2)/64) \sqrt{2\pi \frac{8}{\beta}}\exp(-\beta(\frac{r^2}{2}-2)^2/16)=\frac{\sqrt \pi}{\sqrt \beta}.
\nneq
When $0\le r\le 2$, by (\ref{PE1}) we have
$$4h^{\prime\prime}(r)=-r-\beta r(1-r^2/4)h^{\prime}(r)\le 0,$$
hence
\bbeq\label{cw2}  h^{\prime}(r)\le h^{\prime}(0)=\frac14 e^{\beta/4} \int_{0}^{+\infty}\exp(-\beta(u-2)^2/16)du <  e^{\beta/4}\frac{\sqrt \pi}{\sqrt \beta}.
\nneq
Combining (\ref{cw1}) and (\ref{cw2}), we obtain $\|h^\prime\|_{\infty} < e^{\beta/4}\frac{\sqrt \pi}{\sqrt \beta}$.
Thus assumption $(\bf H)$ holds once if
\bbeq\label{cw3}
|K|\sqrt{\pi \beta} e^{\beta/4} \le 1 \nneq
and then the conclusions of Theorem \ref{W1} and Theorem \ref{W1c} hold under (\ref{cw3}).

For the result of propagation of chaos, we can take $c_1=|K|\beta+\vep^{\prime},\ c_2=\frac{\beta}{4}(1+|K|+\frac{\vep^{\prime}}{\beta})^2$ for any $\vep^{\prime}>0$, and $c_3=0$. Then condition (\ref{MEcc}) is satisfied and then the conclusion of Theorem \ref{Wp} holds under (\ref{cw3}).
}
\nexa

\bexa\label{Double-Well} {\bf (Double-Well confinement potential and quadratic interaction)} Let $d=1$, and $$V(x)=\beta(x^4/4 -x^2/2),\ W(x,y)=\beta K (x-y)^2$$ where $\beta>0$ is the inverse temperature, $K\in \rr$. This model has the double-well confinement potential and quadratic interaction potential.

{\rm
For this model,  we have
$$b_0(r)=\beta r(1-2K-r^2/4), \ \forall r>0.$$
So conditions (\ref{b1}) and (\ref{b2}) are satisfied.
By the similar calculations as in Example \ref{Curie-Weiss}, we get
 $$\|h^\prime\|_{\infty} < \begin{cases}
e^{(1-2K)^2\beta/4}\frac{\sqrt \pi}{\sqrt \beta},\ &\text{ if } K\le \frac12,\\
\frac{\sqrt \pi}{\sqrt \beta},\ &\text{ if } K > \frac12.
\end{cases}$$
Since $\|\nabla_{xy}^2 W\|_{\infty}=2|K|\beta$, assumption $(\bf H)$ holds once if

 \bbeq\label{cw4}\begin{cases}
 2|K| \sqrt{\pi \beta} e^{(1-2K)^2\beta/4} \le 1,  &\text{ if }\ K\le \frac12\\
 2|K| \sqrt{\pi \beta} \le 1, &\text{ if }\ K > \frac12,
 \end{cases}\nneq and then the conclusion of Theorem \ref{W1} holds under (\ref{cw4}).

Furthermore, note that (\ref{b5}) holds with $M=\beta (1-2K)$, and
$$M+\|\nabla_{xy}^2 W\|_{\infty}=\begin{cases}
\beta,\ &\text{ if } K\ge 0,\\
\beta(1-4K),\ &\text{ if } K< 0
\end{cases}$$
which is strictly positive. Then the conclusion of Theorem \ref{W1c} holds.

For the result of propagation of chaos in Theorem \ref{Wp}, we can take $c_3=0$ when $K\ge 0$, and $c_3=-2K\beta$ when $K<0$. To ensure that conditions (\ref{MEaa}) and (\ref{MEcc}) are satisfied, one can take $c_1=2|K|\beta+\vep^{\prime},\ c_2=\frac{\beta}{4}(1+2|K|+\frac{\vep^{\prime}}{\beta})^2$ in the case of $K>0$ and $c_1=-4K\beta+\vep^{\prime},\ c_2=\frac{\beta}{4}(1-4K+\frac{\vep^{\prime}}{\beta})^2$ in the case of $K<0$, for any $\vep^{\prime}>0$.
}\nexa

\section{Proofs of Theorems \ref{W1} and \ref{W1c}}
\subsection{Coupling}
We first introduce the approximate componentwise reflection coupling by following A. Eberle \cite{[EA]}.
Given $\delta>0$, let $\lambda_{\delta}, \pi_{\delta} : \mathbb{R}^{+} \rightarrow [0,1]$ be two Lipschitz continuous functions such that
\begin{equation}\label{CC1} \lambda_{\delta}(r)^2+\pi_{\delta}(r)^2=1,\ \forall r \in \mathbb{R}^{+}
\end{equation}
and
\begin{equation}\label{CC2} \lambda_{\delta}(r)=
\begin{cases} 1,\  \text{if} \ r \ge \delta,\\
0,\  \text{if} \ r \le \delta/2.
\end{cases}
\end{equation}

Then a coupling of two solutions of the mean-field interacting particle system (\ref{MFS}) with initial values $x_0,y_0\in(\rr^d)^N$ is given
by a strong solution of the system
\begin{equation}\label{coupling}
\begin{aligned}
dX_t^{i,N}&=\sqrt{2}[\lambda_{\delta}(|Z_t^i|) dB_t^{1,i} + \pi_{\delta}(|Z_t^i|) dB_t^{2,i}] -\nabla V(X_t^{i,N})dt\\
&\ \ \ \ \ \ \ \  \ \ \ \ \ \ \ \ \ \ \ \  \ \ \ \  \ \ \ \ \ \ \ \ \ \ \ \ \ - \frac 1{N-1} \sum_{j: j\ne i, 1\le j \le N}\nabla_x W(X^{i,N}_t, X^{j,N}_t) dt, \\
dY_t^{i,N}&=\sqrt{2}[\lambda_{\delta}(|Z_t^i|)R_t^i dB_t^{1,i} + \pi_{\delta}(|Z_t^i|) dB_t^{2,i}] -\nabla V(Y_t^{i,N})dt\\
&\ \ \ \ \ \ \ \  \ \ \ \ \ \ \ \ \ \ \ \  \ \ \ \  \ \ \ \ \ \ \ \ \ \ \ \ \ - \frac 1{N-1} \sum_{j: j\ne i, 1\le j \le N}\nabla_x W(Y^{i,N}_t, Y^{j,N}_t) dt, \\
\end{aligned}
\end{equation}
$1\le i\le N.$ Here $Z_t^i:=X_t^{i,N}-Y_t^{i,N}$ and $R_t^i:=I_d-2e_t^i (e_t^{i})^T$, where $I_d$ is the $d$-dimensional unit matrix
and $e_t^i (e_t^{i})^T$ is the orthogonal projection onto the unit vector $e_t^i:=Z_t^i/|Z_t^i|$ if $|Z_t^i|\neq 0$. $B_t^{1,i}$ and $B_t^{2,i}, 1 \le i \le N,$ are independent standard Brownian motions taking values in $\mathbb{R}^{d}$. We will denote $X_t^{(N)}=(X_t^{1,N}, \cdots, X_t^{N,N})$, $Y_t^{(N)}=(Y_t^{1,N}, \cdots, Y_t^{N,N})$ and $Z_t^{(N)}:=X_t^{(N)}-Y_t^{(N)}$.

To see that $(X_t^{(N)},Y_t^{(N)})$ is a coupling process,  it is enough to notice that
\begin{equation}\label{BMs}
\begin{aligned}
\hat{B}_t^i:&=\int_0^t \lambda_{\delta}(|Z_s^i|) dB_s^{1,i} + \int_0^t \pi_{\delta}(|Z_s^i|) dB_s^{2,i}\\
\check{B}_t^i:&=\int_0^t \lambda_{\delta}(|Z_s^i|)R_t^i dB_s^{1,i} + \int_0^t \pi_{\delta}(|Z_s^i|) dB_s^{2,i}, 1\le i \le N,\\
\end{aligned}
\end{equation}
are standard Brownian motions on $(\rr^d)^N$.
\begin{remarks}\label{C5}{\rm  \begin{enumerate}[$(1)$]

\item The coupling (\ref{coupling}) behaves as a reflection coupling when the distance between the two particles $X_t^{i,N}$ and $Y_t^{i,N}$ are larger than $\delta$. When the particles are very close (with distance less than $\frac12 \delta$), they are driven by the same Brownian motion, i.e., it is a synchronous coupling. And when the distance is between $\frac12\delta$ and $\delta$, it is a mixture of reflection coupling and synchronous coupling. The aim is to make $\lambda_{\delta}$ and $\pi_{\delta}$ globally Lipschitz continuous, so that the coupling SDE has a unique strong solution, given the independent Brownian motions $B^{1,i}_t, B^{2,i}_t$, $1\le i\le N$.

 \item If one adopts the componentwise reflection coupling (i.e. the limit coupling when $\delta\to 0$), since $X^{i,N}$, $Y^{i,N}$ will separate after the time that they meet (i.e. $X_t^{i,N}=Y^{i,N}_t$), the local times will appear when It\^o's formula is applied for $|X^{i,N}_t-Y^{i,N}_t|$. This makes the control of $\sum_{i=1}^N |X^{i,N}_t-Y^{i,N}_t|$ difficult to deal with. That is the reason  why A. Eberle \cite{[EA]} introduced the synchronous coupling when $|X^{i,N}_t-Y^{i,N}_t|$ is small.

 \end{enumerate}
}
\nrmk

\subsection{Proofs of Theorem \ref{W1}}

\begin{proof}[Proof of Theorem \ref{W1}] {\bf 1). Proof of (\ref{W2}).} The first inequality in (\ref{W2}) is trivial, and next we prove the second inequality. By doing subtraction of the equations in (\ref{coupling}), we have
\begin{equation}\label{CZ1}
\begin{aligned}
dZ_t^i&=2\sqrt{2}\lambda_{\delta}(|Z_t^i|)e_t^i d\tilde{B}_t^{i}-[\nabla V(X_t^{i,N})-\nabla V(Y_t^{i,N})]dt\\
&\ \ \ \ \ \ \ \ \ \ \ \ -\frac{1}{N-1}\sum_{j: j\ne i, 1\le j \le N}[\nabla_x W(X^{i,N}_t, X^{j,N}_t)-\nabla_x W(Y^{i,N}_t, Y^{j,N}_t)]dt,\\
Z_0^i&=x_0^i-y_0^i,\\
\end{aligned}\end{equation}
where the processes
$\tilde{B}_t^{i}=\int_0^t (e_s^i)^T dB_s^{1,i}, 1\le i\le N,$
are  one-dimensional standard Brownian motions such that $\<\tilde B^i, \tilde B^j\>_t=0$ for $i\ne j$.

Let $r_t^i=|Z_t^i|,  1\le i\le N.$ By applying It\^o's formula, we have

\begin{equation}\label{CZ2}
\begin{aligned}
dr_t^i&=1_{\{r_t^i \neq 0\}}  2\sqrt{2}\lambda_{\delta}(r_t^i) d\tilde{B}_t^{i}-1_{\{r_t^i \neq 0\}} \langle e_t^i, \nabla V(X_t^{i,N})-\nabla V(Y_t^{i,N})\rangle dt \\
&-1_{\{r_t^i \neq 0\}} \langle e_t^i,\frac{1}{N-1}\sum_{j: j\ne i, 1\le j \le N}[\nabla_x W(X^{i,N}_t, X^{j,N}_t)-\nabla_x W(Y^{i,N}_t, Y^{j,N}_t)]\rangle dt\\
&+1_{\{r_t^i \neq 0\}} \sum_{k,l=1}^d [1_{\{k=l\}}(r_t^i)^{-1}
-(X_t^{i,N,k}-Y_t^{i,N,k})(X_t^{i,N,l}-Y_t^{i,N,l})(r_t^i)^{-3}]\lambda_{\delta}(r_t^i)^2(I_d-R_t^i)^2_{kl}dt ,\\
\end{aligned}\end{equation}where $X_t^{i,N,k}$ and $Y_t^{i,N,k}$ denote the $k$-th coordinate of $X_t^{i,N}$ and $Y_t^{i,N}$ respectively, $1\le k\le d$.
Notice that the last term in the right hand side of the above equation equals to 0 by an easy calculation.
Hence we get
\begin{equation}\label{CZ3}
\begin{aligned}
dr_t^i&=1_{\{r_t^i \neq 0\}} 2\sqrt{2}\lambda_{\delta}(r_t^i) d\tilde{B}_t^{i}- 1_{\{r_t^i \neq 0\}}\langle e_t^i,
\frac{1}{N-1}\sum_{j: j\ne i, 1\le j \le N}[\nabla_x W(X^{i,N}_t, X^{j,N}_t)-\nabla_x W(X^{i,N}_t, Y^{j,N}_t)] \rangle dt\\
&-1_{\{r_t^i \neq 0\}}\langle e_t^i, \nabla V(X_t^{i,N})-\nabla V(Y_t^{i,N})+\frac{1}{N-1}\sum_{j: j\ne i, 1\le j \le N}
[\nabla_x W(X^{i,N}_t, Y^{j,N}_t)-\nabla_x W(Y^{i,N}_t, Y^{j,N}_t)]\rangle dt\\
&\le 1_{\{r_t^i \neq 0\}} 2\sqrt{2}\lambda_{\delta}(r_t^i) d\tilde{B}_t^{i}  + \frac{1}{N-1}
\|\nabla^2_{xy}W\|_{\infty}\sum_{j: j\ne i, 1\le j \le N}r_t^j dt + 1_{\{r_t^i \neq 0\}} b_0(r_t^i)dt,\\
\end{aligned}\end{equation}where we use the definition (\ref{b0}) of $b_0$ in the last inequality. Here $d\xi_t\le d\eta_t$ means that $\eta_t-\xi_t$ is a non-decreasing process.

Let $L_{\lambda_{\delta}}$ be the generator defined by for any function $f\in C^2(0,+\infty)$ and $r>0$,
\begin{equation}\label{Ge1}
L_{\lambda_{\delta}}f(r):=4\lambda_{\delta}^2(r)f^{\prime\prime}(r)+b_0(r)f^{\prime}(r).
\end{equation}
Note that $L_{\lambda_{\delta}}$ equals $\LL_{ref}$ when $\lambda_{\delta}\equiv 1$.

Applying It\^o's formula to the function $h(r_t^i)$ and using (\ref{CZ3}) and the fact that $h'(r)>0$, we get for any $t>0$ and $i=1,\cdots, N$,
\begin{equation}\label{h1}
\begin{aligned}
dh(r_t^i)&\le 2\sqrt{2}\lambda_{\delta}(r_t^i) h^{\prime}(r_t^i) d\tilde{B}_t^{i} +h^{\prime}(r_t^i) b_0(r_t^i)dt + 4 h^{\prime\prime} \lambda_{\delta}(r_t^i)^2 dt\\
&\ \ \ \ \ \ \ \ \ \ \ \ \ \ \ \ \ \ \ \ \ \ \ \ \ \ \ \ \ \  \ \ \ \ \ \ \ \ \ \ \ \ + \frac{1}{N-1}\|\nabla^2_{xy}W\|_{\infty}h^{\prime}(r_t^i)\sum_{j: j\ne i, 1\le j \le N}r_t^j dt  \\
&= 2\sqrt{2}\lambda_{\delta}(r_t^i) h^{\prime}(r_t^i) d\tilde{B}_t^{i} + L_{\lambda_{\delta}}h(r_t^i) dt +
\frac{1}{N-1}\|\nabla^2_{xy}W\|_{\infty}h^{\prime}(r_t^i)\sum_{j: j\ne i, 1\le j \le N}r_t^j dt. \\
\end{aligned}\end{equation}

Notice that by the definition of $L_{\lambda_{\delta}}$ and the Poisson equation (\ref{PE1}),
\bbeq\label{Ge2} L_{\lambda_{\delta}}h(r)=\LL_{ref} h(r)+4(\lambda_{\delta}^2 -1)h^{\prime\prime}(r)=-r+(1-\lambda_{\delta}^2)(r+b_0(r)h^{\prime}(r)).\nneq
Then
$$\aligned
-\sum_{i=1}^N &\left( L_{\lambda_{\delta}}h(r_t^i) +
\frac{1}{N-1} \|\nabla^2_{xy}W\|_{\infty} h^{\prime}(r_t^i) \sum_{j: j\ne i, 1\le j \le N}r_t^j \right)\\
&\ge (1-\|\nabla^2_{xy}W\|_{\infty} \|h^{\prime}\|_\infty) \sum_{i=1}^N r_t^i  - \sum_{i=1}^N(1-\lambda_{\delta}(r_t^i)^2)(r_t^i+b_0(r_t^i)h^{\prime}(r_t^i))
\endaligned$$
which is bounded from below by $-N(\delta+ \sup_{r\in (0,\delta)} b_0^+(r)\|h'\|_\infty)$ according to the conditions {\bf (H)} and (\ref{b2}). By integrating from $0$ to $T$ and taking expectation in the previous inequality (\ref{h1}) for $dh(r_t^i)$ and using Fatou's lemma, we have for any $T>0$,

\begin{equation}\label{h3}
\begin{aligned}
&\mathbb{E}\int_0^T \left\{(1-\|\nabla^2_{xy}W\|_{\infty} \|h^{\prime}\|_{\infty})\sum_{i=1}^N r_t^i
-\sum_{i=1}^N(1-\lambda_{\delta}(r_t^i)^2)(r_t^i+b_0(r_t^i)h^{\prime}(r_t^i))\right\}dt \\
&\le \sum_{i=1}^N h(|x_0^i-y_0^i|).\\
\end{aligned}\end{equation}
Letting $\pp_x|_{[0,T]}$ be the law of $(X^{(N)}_t)_{t\in [0,T]}$, we obtain by assumption $(\bf H)$ and (\ref{h3})
\begin{equation}\label{h4}
\begin{aligned}
&W_{1,d_{L^1[0,T]}}(\pp_{x_0}|_{[0,T]}, \pp_{y_0}|_{[0,T]}) \le \mathbb{E}\int_0^T d_{l^1}(X_t^{(N)}, Y_t^{(N)}) dt=\mathbb{E}\int_0^T \sum_{i=1}^N r_t^i dt\\
&\le \frac{1}{1-\|\nabla^2_{xy}W\|_{\infty} \|h^{\prime}\|_{\infty}} \left\{\sum_{i=1}^N h(|x_0^i-y_0^i|)
+ \sum_{i=1}^N \mathbb{E}\int_0^T (1-\lambda_{\delta}(r_t^i)^2)(r_t^i+b_0^{+}(r_t^i)h^{\prime}(r_t^i)) dt \right\}.\\
\end{aligned}\end{equation}

By the definition of $\lambda_{\delta}$ and the assumption $\lim_{r\rightarrow 0}b^+_0(r)=0$,
the second term in the right hand side of the inequality above converges to 0, a.s., as $\delta \downarrow 0$. Hence
\begin{equation}\label{h5}
W_{1,d_{L^1[0,T]}}(\pp_{x_0}|_{[0,T]}, \pp_{y_0}|_{[0,T]}) \le \frac{1}{1-\|\nabla^2_{xy}W\|_{\infty} \|h^{\prime}\|_{\infty}} \sum_{i=1}^N h(|x_0^i-y_0^i|).
\end{equation}

Let $Q_n$ be an optimization coupling of $(\pp_{x_0}|_{[0,n]}, \pp_{y_0}|_{[0,n]})$ for $W_{1,d_{L^1[0,n]}}(\pp_{x_0}|_{[0,n]}, \pp_{y_0}|_{[0,n]})$. Then $\{Q_n|_{[0,T]}; n\ge T\}$ is tight for any finite time $T$ (because their marginal distributions are respectively $\pp_{x_0}|_{[0,T]}$ and $\pp_{y_0}|_{[0,T]}$), hence one can find a probability measure $Q$ on $C(\rr^+, (\rr^d)^N)^2$ such that $Q_n|_{[0,T]}
\to Q|_{[0,T]}$ weakly for all $T>0$. Thus
$$\aligned
W_{1,d_{L^1[0,\infty]}}(\pp_{x_0}, \pp_{y_0})&\le \ee^{Q} \int_0^\infty d_{l^1}(\gamma_1(t), \gamma_2(t)) dt\\
&=\lim_{T\to +\infty} \ee^{Q} \int_0^T d_{l^1}(\gamma_1(t), \gamma_2(t)) dt \\
&\le \lim_{T\to\infty} \lim_{T\le n\to+\infty}\ee^{Q_n} \int_0^T d_{l^1}(\gamma_1(t), \gamma_2(t)) dt\\
&\le \lim_{n\to\infty} W_{1,d_{L^1[0,n]}} (\pp_{x_0}|_{[0,n]}, \pp_{y_0}|_{[0,n]}).
\endaligned$$
The converse inequality is evident. Therefore we have
$$
W_{1,d_{L^1}}(\pp_{x_0}, \pp_{y_0})=\lim_{n\to\infty} W_{1,d_{L^1[0,n]}} (\pp_{x_0}|_{[0,n]}, \pp_{y_0}|_{[0,n]}).
$$
From this and (\ref{h5}) we obtain (\ref{W2}).

{\bf 2). Proof of (\ref{cor1}).} Note that for any Lipschitzian function $g$ w.r.t the $d_{l^1}$-metric on $(\rr^d)^N$, $g$ is $\mu^{(N)}$-integrable because $\int \sum_{i=1}^N |x^i|d\mu^{(N)}(x)<+\infty$. So we can assume $\mu^{(N)}(g)=0$ without loss of generality.  Moreover we have
 $$
 \aligned
 \int_0^{+\infty} | P^{(N)}_tg(x)|dt&=
 \int_0^{+\infty} | P^{(N)}_tg(x)-\int P_t^{(N)}g(y) d\mu^{(N)}(y)| dt\\
 &\le \|g\|_{Lip(d_{l^1})}\int_{(\rr^d)^N}  \int_0^{+\infty} W_{d_{l^1}}(P_t^{(N)}(x,\cdot), P_t^{(N)}(y,\cdot)) dt d\mu^{(N)}(y)\\
&\le \frac{\|g\|_{Lip(d_{l^1})}}{1-\|\nabla^2_{xy}W\|_{\infty} \|h^{\prime}\|_{\infty}} \int_{(\rr^d)^N}  \sum_{i=1}^N h(|x^i-y^i|) d\mu^{(N)}(y)\\
&<+\infty,
 \endaligned$$
 then the unique solution of the Poisson equation $-\LL^{(N)} G=g $ with $\mu^{(N)}(G)=0$ is given by
 $G(x)=\int_0^{+\infty} P^{(N)}_t g(x)dt, \forall x\in(\rr^d)^N$.

 For each $1\le i\le N$, letting $\tilde x^i\ne x^i$ and $\tilde x\in (\rr^d)^N$ so that $(\tilde x)^j=x^j$ for $j\ne i$ and  $(\tilde x)^i=\tilde x^i$, we have
\begin{equation}\label{cor2}
\begin{aligned}
|\nabla_i G(x)|&\le \limsup_{\tilde{x}^i\rightarrow x^i} \frac{|G(x)-G(\tilde{x})|}{|x^i-\tilde{x}^i|}\\
&\le \limsup_{\tilde{x}^i\rightarrow x^i} \frac{1}{|x^i-\tilde{x}^i|} \int_0^{+\infty} | P^{(N)}_tg(x)-P_t^{(N)}g(\tilde{x})| dt\\
&\le \limsup_{\tilde{x}^i\rightarrow x^i} \frac{1}{|x^i-\tilde{x}^i|} \|g\|_{Lip(d_{l^1})} \int_0^{+\infty} W_{d_{l^1}}(P_t^{(N)}(x,\cdot), P_t^{(N)}(\tilde{x},\cdot)) dt\\
&\le \frac{1}{1-\|\nabla^2_{xy}W\|_{\infty} \|h^{'}\|_{\infty}} \|g\|_{Lip(d_{l^1})}
\lim_{\tilde{x}^i\rightarrow x^i}\frac{h(|x^i-\tilde{x}^i|)}{|x^i-y^i|}\\
&= \frac{h^{\prime}(0)}{1-\|\nabla^2_{xy}W\|_{\infty} \|h^{\prime}\|_{\infty}} \|g\|_{Lip(d_{l^1})},  \\
\end{aligned}\end{equation}
where the fourth inequality follows from (\ref{h5}).
\nprf

\subsection{Proof of Theorem \ref{W1c}}
\bprf Here we also adopt the coupling (\ref{coupling}) . Let $h$ be defined as in (\ref{PE2}). Define for any $\vep >0$,
\bbeq\label{h0}h_\vep (r):= h(r)+ \vep r, \forall r \ge 0,\nneq
and
$$H_t^{\vep}:=e^{K_{\vep}t}\sum_{i=1}^N h_{\vep}(r_t^i),$$
where $r_t^i=|X_t^{i,N}-Y_t^{i,N}|, 1\le i\le N$, as in the proof of Theorem \ref{W1}.
By using Ito's formula and (\ref{CZ3}), we get for any $t\ge 0$,
\begin{equation}\label{W1c2}
\begin{aligned}
dH_t^{\vep} &\le 2\sqrt{2}e^{K_{\vep}t}\sum_{i=1}^N  \lambda_{\delta}(r_t^i) d\tilde{B}_t^{i}+K_{\vep}H_t^{\vep}dt
+ e^{K_{\vep}t} \sum_{i=1}^N (L_{\lambda_{\delta}}h(r_t^i)+\vep b_0(r_t^i))dt\\
&\ \ \ \ \ \ \ \ \ \ \ \ \ \ \ \ \ \ \ \ \ \ \ +e^{K_{\vep}t} \sum_{i=1}^N (h^{\prime}(r_t^i)+\vep)\sum_{j: j\ne i}\frac{1}{N-1}\|\nabla^2_{xy}W\|_{\infty}r_t^j dt\\
\end{aligned}\end{equation}

Let
\begin{equation}\label{W1c3}
\begin{aligned}
D_t^{\vep}:&=K_{\vep}H_t^{\vep}+e^{K_{\vep}t} \sum_{i=1}^N (L_{\lambda_{\delta}}h(r_t^i)+\vep b_0(r_t^i))
+\frac{1}{N-1} e^{K_{\vep}t} \|\nabla^2_{xy}W\|_{\infty}   \sum_{ i\ne j, 1\le i, j \le N}(h^{\prime}(r_t^i)+\vep)r_t^j
\end{aligned}\end{equation}
be the drift term at the right hand side above.
Calculating as in the proof of Theorem \ref{W1}, we have
\begin{equation}\label{W1c4}
\begin{aligned}
D_t^{\vep} &\le e^{K_{\vep}t} \sum_{i=1}^N [1-\lambda_{\delta}(r_t^i)^2][r_t^i+b_0(r_t^i)h^{\prime}(r_t^i)]\\
&\ \ \ \ \ \ \ \ \ \ \ \ \ \ +e^{K_{\vep}t} \sum_{i=1}^N \{K_{\vep}h_{\vep}(r_t^i)-[1-(\|h^{\prime}\|_{\infty}+\vep)\|\nabla^2_{xy}W\|_{\infty}]r_t^i
+ \vep b_0(r_t^i)\}\\
&\le e^{K_{\vep}t} \sum_{i=1}^N [1-\lambda_{\delta}(r_t^i)^2][r_t^i+b_0(r_t^i)h^{\prime}(r_t^i)]\\
&\ \ \ \ \ \ \ \ \ \ \ \ \ \ +e^{K_{\vep}t} \sum_{i=1}^N \{K_{\vep}(\|h^{\prime}\|_{\infty}+\vep)+\vep M
-[1-(\|h^{\prime}\|_{\infty}+\vep)\|\nabla^2_{xy}W\|_{\infty}] \}r_t^i, \\
\end{aligned}\end{equation}
where we use the assumption $b_0(r)\le Mr,\forall r>0$.

By taking \bbeq \label{W1ca}K_{\vep}=\frac{1-\|\nabla^2_{xy}W\|_{\infty}\|h^{\prime}\|_{\infty}- \vep (M+\|\nabla^2_{xy}W\|_{\infty})}{\|h^{\prime}\|_{\infty}+\vep},\nneq the second term in the right hand of the inequality above vanishes. Then by taking expectation in (\ref{W1c2}) and the localization stopping time technique, we have for any $t\ge 0$,
\begin{equation}\label{W1c5}
\begin{aligned}
\mathbb{E} e^{K_{\vep}t}\sum_{i=1}^N h_{\vep}(r_t^i) \le \sum_{i=1}^N h_{\vep}(|x_0^i-y_0^i|)
+ \mathbb{E}\int_0^t e^{K_{\vep}s}[1-\lambda_{\delta}(r_t^i)^2][r_t^i+b_0^+(r_t^i)h^{\prime}(r_t^i)] ds .\\
\end{aligned}\end{equation}

Note that the second term in the right hand side of the above inequality converges to $0$ as $\delta\downarrow 0$, under the assumption (\ref{b2}). Therefore we get
\begin{equation}\label{W1c6}
\aligned
W_{1,d_{l^1}}(P^{(N)}_t(x_0,\cdot), & P^{(N)}_t(y_0,\cdot)) \le \lim_{\delta\to0}\ee \sum_{i=1}^N r_t^i\\
&\le \sup_{r>0}\frac{r}{h(r)+\vep r}  \lim_{\delta\to0} \mathbb{E} \sum_{i=1}^N h_{\vep}(r_t^i) \\
&\le \sup_{r>0}\frac{r}{h(r)+\vep r}  e^{-K_{\vep}t}\sum_{i=1}^N h_{\vep}(|x_0^i-y_0^i|)\\
&\le \sup_{r>0}\frac{r}{h(r)+\vep r} \cdot \sup_{r>0}\frac{h(r)+\vep r}{r}  e^{-K_{\vep}t}\sum_{i=1}^N|x_0^i-y_0^i|
\endaligned \end{equation}
where the third inequality above follows by (\ref{W1c5}).
\nprf

\section{Propagation of chaos}
We begin with a uniform in time control of the second moment, which is more or less known, see e.g. Cattiaux {\it et al.} \cite{[CGM]}.
\blem\label{ME} Suppose that there exist some positive constants $c_1,c_2,c_3$ such that
 \bbeq\label{MEa} \langle x, \nabla V(x)\rangle\ge c_1|x|^2-c_2, \forall x\in\rr^d
 \nneq
 and
 \bbeq\label{MEb} \langle z, \nabla^2_{xx}W(x,y) z\rangle\ge -c_3 |z|^2, \forall x,y,z\in\rr^d.
 \nneq
 Assume
\bbeq\label{MEc}c_1-c_3-\|\nabla^2_{xy}W\|_{\infty}>0.\nneq
 Let $X_t$ be a solution of (\ref{SDE}) with $\mathbb{E}|X_0|^2<\infty$, then for any
 $\vep\in (0, c_1-c_3-\|\nabla^2_{xy}W\|_{\infty})$,
 \bbeq\label{MEd}\sup_{t\ge0}\mathbb{E}(|X_t|^2)^{\frac 12}\le \max\{ m_2(\mu_0), \hat c(\vep)\},\nneq
 where $m_2(\mu_0)$ and $\hat c(\vep)$ are given in \eqref{mc}.
\nlem
\bprf By It\^o's formula, we have
$$ \begin{aligned} d|X_t|^2&=-2\langle X_t, \nabla V(X_t) \rangle dt  -2\langle X_t, \nabla_x W\circledast\mu_t(X_t) \rangle dt +2d \cdot dt+ 2\sqrt{2}\<X_t, dB_t\>\\
\end{aligned}$$

Notice that for any $x\in \rr^d$, we have
\bbeq\label{ME2} \begin{aligned}
\langle x, \nabla_x W\circledast\mu_t(x)-\nabla_x W\circledast\mu_t(0)\rangle &=\langle x,\int_0^1 \frac{d}{ds} \nabla_x W\circledast\mu_t(sx) ds \rangle\\
&=\langle x, \int_0^1 \frac{d}{ds} \int_{\rr^d}\nabla_x W(sx,y)\mu_t(dy) ds \rangle\\
&=\int_0^1 \int_{\rr^d} \langle x, \nabla_{xx}^2 W(sx,y) x\rangle \mu_t(dy) ds\\
&\ge -c_3 |x|^2,\\
\end{aligned}\nneq where the last inequality follows from (\ref{MEb}).

On the other hand,
\bbeq\label{ME3} \begin{aligned}
|\nabla_x W\circledast\mu_t(0)\rangle |&\le |\nabla_x W(0,0)| +\int_{\rr^d}|\nabla_x W(0,y)-\nabla_x W(0,0)|\mu_t(dy)\\
&\le |\nabla_x W(0,0)| + \|\nabla^2_{xy}W\|_{\infty} \ee |X_t|.\\
\end{aligned}\nneq
Therefore we have
$$\begin{aligned} d |X_t|^2&\le2\left(c_3 |X_t|^2 +\|\nabla^2_{xy}W\|_{\infty}  |X_t| \ee|X_t| + |\nabla_x W(0,0)| |X_t| \right)dt\\
&\quad \quad \quad +2( -c_1 |X_t|^2 + c_2+d)dt+2\sqrt{2}\<X_t, dB_t\>\\
&\le -2(c_1-c_3-\vep) |X_t|^2dt +2\|\nabla^2_{xy}W\|_{\infty}  |X_t| \ee|X_t|dt\\
&\quad \quad \quad + 2(d+c_2+ \frac 1{4\vep} |\nabla_x W(0,0)|^2)dt+2\sqrt{2}\<X_t, dB_t\>\\
\end{aligned}$$
where $0<\vep < c_1-c_3-\|\nabla^2_{xy}W\|_{\infty}$. By the previous stochastic differential  inequality,
$$
|X_t|^2+\int_0^t [2(c_1-c_3-\vep) |X_s|^2 -2\|\nabla^2_{xy}W\|_{\infty}  |X_s| \ee|X_s|]ds- 2t(d+c_2+ \frac 1{4\vep} |\nabla_x W(0,0)|^2)
$$
is a local supermartingale, then a supermartingale by Fatou's lemma. Then for any $T>0$, we have
 $$\begin{aligned} \ee |X_0|^2&\ge \ee |X_T|^2+ 2(c_1-c_3-\vep)\int_0^T  \ee|X_s|^2 ds-2\|\nabla^2_{xy}W\|_{\infty} \int_0^T (\ee|X_s|)^2 ds\\
 &\quad\quad\quad -2T(d+c_2+ \frac 1{4\vep} |\nabla_x W(0,0)|^2)\\
 &\ge 2(c_1-c_3-\vep-\|\nabla^2_{xy}W\|_{\infty})\int_0^T  \ee|X_s|^2 ds-2T(d+c_2+ \frac 1{4\vep} |\nabla_x W(0,0)|^2),\\
 \end{aligned}$$
 which implies $\ee \int_0^T |X_s|^2 ds<+\infty$. In other words
$\int_0^t 2\sqrt{2}\<X_s, dB_s\>$ is a $L^2$-martingale. By taking expectation in (\ref{MEa}) we obtain by (\ref{ME2}) and (\ref{ME3}),
\bbeq\label{ME1} \begin{aligned} \frac {d}{dt} \ee |X_t|^2&\le -2c_1\ee |X_t|^2 +2[c_3\ee |X_t|^2 +\|\nabla^2_{xy}W\|_{\infty} (\ee |X_t|)^2 + |\nabla_x W(0,0)| \ee |X_t| ] +2(d+c_2) \\
&\le -2(c_1-c_3-\|\nabla^2_{xy}W\|_{\infty}-\vep)\mathbb{E}|X_t|^2 + 2(d+c_2+ \frac 1{4\vep} |\nabla_x W(0,0)|^2 )
\end{aligned}\nneq
where $0<\vep < c_1-c_3-\|\nabla^2_{xy}W\|_{\infty}$. By Gronwall's lemma we get for any $t \ge 0$
$$\aligned
\ee |X_t|^2 &\le e^{-2(c_1-c_3-\|\nabla^2_{xy}W\|_{\infty}-\vep) t } \left(\ee |X_0|^2 -\frac{d+c_2+ \frac 1{4\vep}|\nabla_x W(0,0)|^2}{c_1-c_3-\|\nabla^2_{xy}W\|_{\infty}-\vep}\right)\\
&\ \ \ \ \ \ \ \ \ + \frac{d+c_2+ \frac 1{4\vep}|\nabla_x W(0,0)|^2}{c_1-c_3-\|\nabla^2_{xy}W\|_{\infty}-\vep}\\
&\le  \max\left\{\ee |X_0|^2, \frac{d+c_2+ \frac 1{4\vep}|\nabla_x W(0,0)|^2}{c_1-c_3-\|\nabla^2_{xy}W\|_{\infty}-\vep}\right\}
\endaligned$$
the desired result.
\nprf

Following the proof above we have the much stronger uniform Gaussian integrability for $X_t$, which should be of independent interest.

\blem\label{lem42} Assume (\ref{MEa}), (\ref{MEb}) and (\ref{MEc}).
 Let $X_t$ be a solution of (\ref{SDE}) with $$\mathbb{E}\exp\left(\lambda_0|X_0|^2\right)<\infty,\quad \text{for some } \lambda_0>0.$$ If
 $$
 0<\lambda\le \min\{\lambda_0; \frac 12 (c_1-c_3-\|\nabla^2_{xy}W\|_{\infty}-\vep) \}$$
 for some $\vep>0$, then
 $$
 \sup_{t\ge0} \ee \exp(\lambda |X_t|^2)<+\infty.
 $$
\nlem
\bprf By It\^o's formula, we have by the estimates leading to (\ref{ME1}) in the proof of Lemma \ref{ME},
$$ \begin{aligned} &d\exp(\lambda |X_t|^2)\\
=&\lambda \exp(\lambda |X_t|^2) \left([2d-2\langle X_t, \nabla V(X_t)+\nabla_x W\circledast\mu_t(X_t) \rangle] dt+
2\sqrt{2}\<X_t, dB_t\> \right)\\
&\ + 4 \lambda^2|X_t|^2 \exp(\lambda |X_t|^2)dt \\
\le & \lambda\exp(\lambda |X_t|^2)\left[ -2(c_1-c_3-\|\nabla^2_{xy}W\|_{\infty}-\vep-2\lambda) |X_t|^2 + 2(d+c_2+ \frac 1{4\vep} |\nabla_x W(0,0)| )\right]dt\\
&+ \lambda \exp(\lambda |X_t|^2)2\sqrt{2}\<X_t, dB_t\>
\end{aligned}$$
where $\vep>0,\lambda>0$ verify $c_1-c_3-\|\nabla^2_{xy}W\|_{\infty}-\vep-2\lambda>0$. Taking $L>0$ large sufficient so that
$$c_5:=2(c_1-c_3-\|\nabla^2_{xy}W\|_{\infty}-\vep-2\lambda) L^2 - 2(d+c_2+ \frac 1{4\vep} |\nabla_x W(0,0)|)>0,
$$
and noting that
$$
-a x^2 +b \le -(a L^2-b) + aL^21_{|x|\le L},\ \forall a>0, \forall x\in\rr,
$$
we obtain by following the same argument as in Lemma \ref{ME}
$$
\frac d{dt}\ee \exp(\lambda |X_t|^2)\le -\lambda c_5\ee \exp(\lambda |X_t|^2) +2(c_1-c_3-\|\nabla^2_{xy}W\|_{\infty}-\vep-2\lambda) L^2 \lambda e^{\lambda L^2}.
$$
Therefore by Gronwall's lemma
$$
\sup_{t\ge0} \ee \exp(\lambda |X_t|^2)<+\infty.
$$
\nprf

Next we present the proof of Theorem \ref{Wp}, which is quite close to those of Theorems \ref{W1} and \ref{W1c} .
\bprf[Proof of Theorem \ref{Wp}]
Let $\lambda_{\delta}$ and $\pi_{\delta}$ be defined as in Section 3.1. Consider the following coupling between the independent copies $\bar{X}_t^i, 1\le i\le N$ of the nonlinear diffusion processes (\ref{SDE}) and the mean-field interacting particle system (\ref{MFS}):
\begin{equation}\label{coupling2}
\begin{aligned}
d\bar{X}_t^{i}&=\sqrt{2}[\lambda_{\delta}(|Z_t^i|) dB_t^{1,i} + \pi_{\delta}(|Z_t^i|) dB_t^{2,i}]-\nabla V(\bar{X}_t^{i})dt - \nabla_x W\circledast \mu_t(\bar{X}_t^{i})dt,\\
dX_t^{i,N}&=\sqrt{2}[\lambda_{\delta}(|Z_t^i|)R_t^i dB_t^{1,i} + \pi_{\delta}(|Z_t^i|) dB_t^{2,i}] -\nabla V(X_t^{i,N})dt\\
&\ \ \ \ \ \ \ \  \ \ \ \ \ \ \ \ \ \ \ \  \ \ \ \  \ \ \ \ \ \ \ \ \ \ \ \ \ - \frac 1{N-1} \sum_{j: j\ne i, 1\le j \le N}\nabla_x W(X^{i,N}_t, X^{j,N}_t) dt. \\
\end{aligned}
\end{equation} Here $Z_t^i:=\bar{X}_t^i-X_t^{i,N}$ and $R_t^i:=I_d-2e_t^i e_t^{i,T}$, where $I_d$ is the $d$-dimensional unit matrix
and $e_t^i e_t^{i,T}$ is the orthogonal projection onto the unit vector $e_t^i:=Z_t^i/|Z_t^i|$ if $|Z_t^i|\neq 0$. $B_t^{1,i}$ and $B_t^{2,i}, 1 \le i \le N,$ are independent standard Brownian motions in $\mathbb{R}^{d}$. We assume that $\bar{X}_t^i$ and $X_t^{i,N}$, $1\le i\le N$ have the same starting points $X_0^{i}, 1\le i\le N$, i.i.d. of law $\mu_0$. The independence of $\bar X^i_t, 1\le i\le N$ comes from the fact that the Brownian motions $\{\int_0^t \lambda_{\delta}(|Z_s^i|) dB_s^{1,i} + \int_0^t\pi_{\delta}(|Z_s^i|) dB_s^{2,i}, \ 1\le i\le N\}$ are independent because their inter-brackets are zero.

By doing subtraction of the equations in (\ref{coupling2}), we have
$$
\begin{aligned}
dZ_t^i&=2\sqrt{2}\lambda_{\delta}(|Z_t^i|)e_t^i d\tilde{B}_t^{i}-[\nabla V(\bar{X}_t^i)-\nabla V(X_t^{i,N})]dt- \nabla_x W\circledast \mu_t(\bar{X}_t^i)dt\\
&\ \ \ \ \ \ \ \ \ \ \ \ \ \ \ \ \ \  \ \ \ \  \ \ \ \ \ \ \ \ \ \ \ \ +\frac 1{N-1} \sum_{j: j\ne i, 1\le j \le N}\nabla_x W(X^{i,N}_t, X^{j,N}_t) dt,\\
\end{aligned}$$
where the processes
$\tilde{B}_t^{i}=\int_0^t (e_s^i)^T dB_s^{1,i}, 1\le i\le N,$
are  one-dimensional standard Brownian motions such that $\<\tilde B^i, \tilde B^j\>_t=0$ for $i\ne j$.

Let $r_t^i=|Z_t^i|,  1\le i\le N.$ By applying It\^o's formula, we have
$$
\begin{aligned}
dr_t^i&=1_{\{r_t^i \neq 0\}}  2\sqrt{2}\lambda_{\delta}(r_t^i) d\tilde{B}_t^{i}-1_{\{r_t^i \neq 0\}} \langle e_t^i, \nabla V(\bar{X}_t^i)
-\nabla V(X_t^{i,N})\rangle dt \\
&\ \ \ \ \ \ \ \ \ \ -1_{\{r_t^i \neq 0\}} \langle e_t^i, \nabla_x W\circledast \mu_t(\bar{X}_t^i)-\frac{1}{N-1}\sum_{j: j\ne i, 1\le j \le N}\nabla_x W(X^{i,N}_t, X^{j,N}_t)\rangle dt\\
&=1_{\{r_t^i \neq 0\}}  2\sqrt{2}\lambda_{\delta}(r_t^i) d\tilde{B}_t^{i}\\
&\ \ \ \ \ \ \ \ \ \ -1_{\{r_t^i \neq 0\}} \langle e_t^i, \nabla V(\bar{X}_t^i) -\nabla V(X_t^{i,N})\rangle dt\\
&\ \ \ \ \ \ \ \ \ \ -1_{\{r_t^i \neq 0\}} \langle e_t^i, \nabla_x W\circledast \mu_t(\bar{X}_t^i)- \frac{1}{N-1}\sum_{j: j\ne i, 1\le j \le N}\nabla_x W(\bar{X}_t^i, \bar{X}_t^j)\rangle dt\\
&\ \ \ \ \ \ \ \ \ \ -1_{\{r_t^i \neq 0\}} \langle e_t^i, \frac{1}{N-1}\sum_{j: j\ne i, 1\le j \le N}[\nabla_x W(\bar{X}_t^i, \bar{X}_t^j)-\nabla_x W(\bar{X}_t^i, X^{j,N}_t)]  \rangle dt\\
&\ \ \ \ \ \ \ \ \ \ -1_{\{r_t^i \neq 0\}} \langle e_t^i, \frac{1}{N-1}\sum_{j: j\ne i, 1\le j \le N}[\nabla_x W(\bar{X}_t^i, X^{j,N}_t)-\nabla_x W(X^{i,N}_t, X^{j,N}_t)] \rangle dt.
\end{aligned}
$$
Remark that the sum of the first and the fourth drift terms above is $\le b_0(r_t^i) dt$, the third drift term above is $\le \frac{1}{N-1}\|\nabla^2_{xy}W\|_{\infty}\sum_{j: j\ne i, 1\le j \le N}r_t^j dt$, and the second drift term is bounded by $I_t^i dt$, where
\bbeq\label{CP1}
I_t^i:=|\nabla_x W\circledast \mu_t(\bar{X}_t^i)- \frac{1}{N-1}\sum_{j: j\ne i, 1\le j \le N}\nabla_x W(\bar{X}_t^i, \bar{X}_t^j)|.
\nneq
Therefore we obtain
\begin{equation}\label{CP2}
dr_t^i \le 2\sqrt{2}\lambda_{\delta}(r_t^i) d\tilde{B}_t^{i} + b_0(r_t^i)dt + \frac{1}{N-1}\|\nabla^2_{xy}W\|_{\infty}\sum_{j: j\ne i, 1\le j \le N}r_t^j dt + I_t^i dt.
\end{equation}

Recall that for any $\vep\ge 0$, $h_{\vep}(r)=h(r)+\vep r, \forall r\ge0$. By using (\ref{CP2}) and It\^o's formula again, we get

\begin{equation}\label{CP3}
\begin{aligned}
dh_{\vep}(r_t^i)&\le 2\sqrt{2}\lambda_{\delta}(r_t^i) h^{\prime}_{\vep}(r_t^i) d\tilde{B}_t^{i}+ 4\lambda_{\delta}^2(r_t^i)h^{\prime\prime}_{\vep}(r_t^i)dt+ b_0(r_t^i)h^{\prime}_{\vep}(r_t^i)dt \\
&\ \ \ \ \ \ \ \ \ \ \ \ \ \ \ \ \ \ \ \ \ \ \ \ \ \ \ \ + \frac{1}{N-1}\|\nabla^2_{xy}W\|_{\infty}h^{\prime}_{\vep}(r_t^i)\sum_{j: j\ne i, 1\le j \le N}r_t^j dt + h^{\prime}_{\vep}(r_t^i)I_t^i dt \\
&=2\sqrt{2}\lambda_{\delta}(r_t^i) h^{\prime}_{\vep}(r_t^i) d\tilde{B}_t^{i}+ [4\lambda_{\delta}^2(r_t^i)h^{\prime\prime}(r_t^i)+b_0(r_t^i)h^{\prime}(r_t^i)]dt+\vep b_0(r_t^i)dt  \\
&\ \ \ \ \ \ \ \ \ \ \ \ \ +\frac{1}{N-1}\|\nabla^2_{xy}W\|_{\infty}(h^{\prime}(r_t^i)+\vep)\sum_{j: j\ne i, 1\le j \le N}r_t^j dt + (h^{\prime}(r_t^i)+\vep)I_t^i dt\\
&\le2\sqrt{2}\lambda_{\delta}(r_t^i) h^{\prime}_{\vep}(r_t^i) d\tilde{B}_t^{i}+ [1-\lambda_{\delta}^2(r_t^i)][r_t^i+b_0(r_t^i)h^{\prime}(r_t^i)]dt -(1- \vep M)r_t^i dt\\
&\ \ \ \ \ \ \ \ \ \ \ \ \ +\frac{1}{N-1}\|\nabla^2_{xy}W\|_{\infty}(\|h^{\prime}\|_{\infty}+\vep)\sum_{j: j\ne i, 1\le j \le N}r_t^j dt + (\|h^{\prime}\|_{\infty}+\vep)I_t^i dt,\\
\end{aligned}\end{equation}
where the last inequality follows from (\ref{Ge1}), (\ref{Ge2}) and (\ref{b5}).

Taking expectation in the inequality above and using the fact that $r_t^i, 1\le i\le N$ have the same law, and setting
$$c_{\vep}:=1-\|\nabla^2_{xy}W\|_{\infty}\|h^{\prime}\|_{\infty}- \vep (M+\|\nabla^2_{xy}W\|_{\infty}),$$
we have

\begin{equation}\label{CP4}
\begin{aligned}
d& \mathbb{E}h_{\vep}(r_t^1)
 \le \mathbb{E}[1-\lambda_{\delta}(r_t^1)^2][r_t^1+b^+_0(r_t^1)h^{\prime}(r_t^1)]dt+(\|h^{\prime}\|_{\infty}+\vep) \mathbb{E}I_t^1 dt-c_{\vep} \mathbb{E} r_t^1dt
\end{aligned}\end{equation}
{\it Proof of part (a).} Choose $\vep=0$, $c_0=1-\|\nabla^2_{xy}W\|_{\infty}\|h^{\prime}\|_{\infty}$. For any $1\le k\le N$, by (\ref{CP4}) we have
$$
\aligned
\frac 1k W_{1,d_{L^1[0,T]}}&(\pp_{\mu_0^{\otimes N}}^{[1,k],N}|_{[0,T]}, \qq_{\mu_0}^{\otimes k}|_{[0,T]})\le \frac 1k \ee \int_0^T\sum_{i=1}^k r_t^i dt = \int_0^T \ee r_t^1dt\\
& \le\frac 1{c_0} \|h^{\prime}\|_{\infty}  \int_0^T \ee I_t^1 dt dt + \frac 1{c_0} \ee \int_0^T [1-\lambda_{\delta}(r_t^1)^2][r_t^1+b^+_0(r_t^1)h^{\prime}(r_t^1)]dt.
\endaligned$$
Letting $\delta\to0+$, the last term tends to zero. Hence
\bbeq\label{CP44}
\frac 1k W_{1,d_{L^1[0,T]}}(\pp_{\mu_0^{\otimes N}}^{[1,k],N}|_{[0,T]}, \qq_{\mu_0}^{\otimes k}|_{[0,T]})\le \frac {\|h^{\prime}\|_{\infty}} {c_0}  \int_0^T \ee I_t^1 dt.
\nneq

Next we estimate $\mathbb{E}I_t^1$, which is the only new point w.r.t. the proofs in Theorems \ref{W1} and \ref{W1c}. Note that $\bar{X}_t^{j}, 2\le j\le N$ are independent copies of $\bar{X}_t^{1}$, and $$\mathbb{E}[\nabla_x W(\bar{X}_t^{1},\bar{X}_t^{j})|\bar{X}_t^{1}]=\nabla_x W\circledast \mu_t(\bar{X}_t^{1}).$$
Thus by using Cauchy-Schwartz inequality, we get
\begin{equation}\label{CP5}
\begin{aligned}
\mathbb{E}I_t^1&\le \left(\mathbb{E}\left\{\mathbb{E}\left[|\nabla_x W\circledast \mu_t(\bar{X}_t^{1})- \frac{1}{N-1}\sum_{2\le j \le N}\nabla_x W(\bar{X}_t^{1}, \bar{X}_t^{j})|^2|\bar{X}_t^{1}\right]\right\}\right)^{\frac12}\\
&=\left(\mathbb{E} \frac{1}{N-1}\int |\nabla_x W(\bar{X}_t^{1},y)- \nabla_x W*\mu_t(\bar{X}_t^{1})|^2 d\mu_t(y) \right)^{\frac12}\\
&\le \frac{1}{\sqrt{N-1}}\|\nabla^2_{xy}W\|_{\infty}\left(\int_{x\in\rr^d}|x-\mu_t(\bar{X}_t^{1})|^2 \mu_t(dx)\right)^{\frac12} \\
&\le \frac{1}{\sqrt{N-1}}\|\nabla^2_{xy}W\|_{\infty}\sup_{t\ge0}(\mathbb{E}|X_t|^2)^{\frac12}.\\
\end{aligned}\end{equation}

Plugging \eqref{CP5} into \eqref{CP44}, we get
$$ \frac 1k W_{1,d_{L^1[0,T]}}(\pp_{\mu_0^{\otimes N}}^{[1,k],N}|_{[0,T]}, \qq_{\mu_0}^{\otimes k}|_{[0,T]})\le  \frac{T}{\sqrt{N-1}}\frac{\|\nabla^2_{xy}W\|_{\infty}\|h^\prime\|_\infty }{1-\|\nabla^2_{xy}W\|_{\infty}\|h^{\prime}\|_{\infty}} \sup_{t\ge0}(\mathbb{E}|X_t|^2)^{\frac12}.
$$
Then by using Lemma \ref{ME}, we obtain the desired result \eqref{MEdd}.

\bigskip
{\it Proof of part (b).} For any $\vep>0$, by \eqref{CP4} we have
\begin{equation}\label{CP4b}
d \mathbb{E}h_{\vep}(r_t^1)
 \le \mathbb{E}[1-\lambda_{\delta}(r_t^1)^2][r_t^1+b^+_0(r_t^1)h^{\prime}(r_t^1)]dt+(\|h^{\prime}\|_{\infty}+\vep) \mathbb{E}I_t^1 dt-c_{\vep} \cdot \inf_{r>0}\frac{r}{h(r)+\vep r} \mathbb{E} h_\vep(r_t^1)dt
\end{equation}

Plugging (\ref{CP5}) into (\ref{CP4b}), we obtain by Gronwall's inequality that for any $\vep>0$ so that $\beta=c_\vep \cdot \inf_{r>0}\frac{r}{h(r)+\vep r}>0$ (i.e. $K_\vep>0$),
\begin{equation}\label{CP6}\aligned
\inf_{r>0} \frac{h_\vep(r)}{ r}\cdot & \ee |\bar X^1_t-X_t^{1,N}|)\le
\mathbb{E}h_{\vep}(|\bar X^1_t-X_t^{1,N}|)\\
\le &\int_0^te^{-\beta  (t-s)}   \frac{1}{\sqrt{N-1}}(\|h^{\prime}\|_{\infty}+\vep)
 \|\nabla^2_{xy}W\|_{\infty}\sup_{t\ge0}(\mathbb{E}|X_t|^2)^{\frac12}ds\\
  &\ + \int_0^te^{-\beta  (t-s)} \mathbb{E}[1-\lambda_{\delta}(r_s^1)^2][r_s^1+b^+_0(r_s^1)h^{\prime}(r_s^1)]ds.\\
\endaligned\end{equation}
By letting $\delta\to0+$, the last term tends to zero. As the joint law of $(\bar X^i_t, 1\le i\le k)$ is $\mu_t^{\otimes k}$, we get for any $1\le k\le N$,
\begin{equation}\label{CP7}\aligned
W_{1, d_{l^1}}(\mu_t^{\otimes k}, \mu_t^{[1,k],N}) \le & \limsup_{\delta\to 0} \ee \sum_{i=1}^k|\bar X^i_t-X_t^{i,N}|= k\cdot\limsup_{\delta\to 0} \ee |\bar X^1_t-X_t^{1,N}| \\
\le &k\cdot \sup_{r>0}\frac{r}{h(r)+\vep r} \frac{1}{\beta \sqrt{N-1}}(\|h^{\prime}\|_{\infty}+\vep) \|\nabla^2_{xy}W\|_{\infty}\sup_{t\ge0}(\mathbb{E}|X_t|^2)^{\frac12}\\
= & \frac{k}{\sqrt{N-1}}\frac{A_\vep}{K_\vep} \|\nabla^2_{xy}W\|_{\infty}\sup_{t\ge0}(\mathbb{E}|X_t|^2)^{\frac12},
\endaligned\end{equation}
which completes the proof by using Lemma \ref{ME}.
\nprf

The proof above yields

\bprop\label{prop43} Under the conditions of Theorem \ref{Wp}, we have
\bbeq\label{prop43a}
\ee W_1\left(\frac 1N \sum_{i=1}^N \delta_{X^{i,N}_t}, \frac 1N \sum_{i=1}^N \delta_{\bar X^{i}_t}\right) \le \frac{1}{\sqrt{N-1}}\frac{A_\vep}{K_\vep} \|\nabla^2_{xy}W\|_{\infty}\sup_{t\ge0}(\mathbb{E}|X_t|^2)^{\frac12}
\nneq
where $(\bar X^i_t)_{t\ge0}, i\ge 1$ are independent copies of the solution $(X_t)_{t\ge0}$ of the McKean-Vlasov equation (\ref{SDE}), and $X_t^{i,N}, 1\le i\le N$  are defined as in \eqref{MFS}.
\nprop

\bprf
Notice that
$$
\ee W_1\left(\frac 1N \sum_{i=1}^N \delta_{X^{i,N}_t}, \frac 1N \sum_{i=1}^N \delta_{\bar X^{i}_t}\right) \le \ee[\frac 1N \sum_{i=1}^N |X^{i,N}_t-\bar X^{i}_t|]= \ee\frac 1N \sum_{i=1}^N r_t^i,
$$ where $r_t^i,1\le i\le N$ are the same as defined in the proof of Theorem \ref{Wp}.
And by (\ref{CP7}), we have
$$
\limsup_{\delta\to 0}\ee \frac 1N\sum_{i=1}^n  r_t^i=\limsup_{\delta\to 0}\ee r_t^1\le \frac{1}{\sqrt{N-1}}\frac{A_\vep}{K_\vep} \|\nabla^2_{xy}W\|_{\infty}\sup_{t\ge0}(\mathbb{E}|X_t|^2)^{\frac12}.
$$
Therefore we obtain (\ref{prop43a}).
\nprf

\brmk{\rm In one-dimensional case, i.e. $d=1$, it is well known that
$$
W_1\left(\frac 1N \sum_{i=1}^N \delta_{ \bar X^{i}_t}, \mu_t\right)=\int_{-\infty}^\infty |\frac 1N \sum_{i=1}^N 1_{(-\infty,x]}( \bar X^i_t)-\mu_t(-\infty, x]|dx.
$$
Then letting $F_t(x)=\mu_t(-\infty, x]$ (the cumulative distribution function), we have by the Cauchy-Schwarz inequality,
$$\aligned
\ee W_1\left(\frac 1N \sum_{i=1}^N \delta_{ \bar X^{i}_t}, \mu_t\right)&\le \int_{\rr} \sqrt{{\rm Var}\left(\frac 1N \sum_{i=1}^N 1_{(-\infty,x]}( \bar X^i_t)\right)}dx\\
&= \frac 1{\sqrt{N}} \int_\rr \sqrt{ F_t(x)(1-F_t(x))} dx,
\endaligned
$$
where the last factor is uniformly bounded in time $t>0$ by some constant $K$ once if $\sup_{t\ge0} \ee |X_t|^{2+\vep}<+\infty$ for some $\vep>0$. The latter uniform $2+\vep$-moment condition is verified once if $\mu_0$ has the $2+\vep$-moment by the arguments in Lemmas \ref{ME} and \ref{lem42}. In other words if $\mu_0$ has the $2+\vep$-moment, there is some constant $K>0$ such that
\bbeq\label{prop43b}
\ee W_1\left(\frac 1N \sum_{i=1}^N \delta_{\bar X^{i}_t}, \mu_t\right)\le \frac K{\sqrt{N}},\ \forall t>0
\nneq
and then the same type bound holds for $\ee W_1\left(\frac 1N \sum_{i=1}^N \delta_{X^{i,N}_t}, \mu_t\right)$, by Proposition \ref{prop43} and the triangular inequality.

But (\ref{prop43b}) does not hold in the multi-dimensional ($d>1$) case, see Fournier and Guillin \cite{[FG]}.
}\nrmk

\brmk\label{Prop43Rem2}{\rm A consequence of Proposition \ref{prop43} is on the bias of $\frac 1N \sum_{i=1}^N f(X^{i,N}_t)$ from $\mu_t(f)$: if $f$ is Lipschitzian on $\rr^d$,
$$\aligned
bias_t(f)&:=|\ee \frac 1N \sum_{i=1}^N f(X^{i,N}_t)- \mu_t(f)|=|\ee \frac 1N \sum_{i=1}^N f(X^{i,N}_t)- \ee \frac 1N \sum_{i=1}^N f(\bar X^{i}_t)| \\
&\le \|f\|_{Lip} \ee W_1 \left(\frac 1N \sum_{i=1}^N \delta_{X^{i,N}_t}, \frac 1N \sum_{i=1}^N \delta_{\bar X^{i}_t}\right)\\
&\le \frac{\|f\|_{Lip}}{\sqrt{N-1}}\frac{A_\vep}{K_\vep} \|\nabla^2_{xy}W\|_{\infty}\sup_{t\ge0}(\mathbb{E}|X_t|^2)^{\frac12}.\endaligned
$$
It is expected that the bias is of order $O(1/N)$, which remains an open question.
}
\nrmk

\section{Quantitative Concentration inequalities}
This section is devoted to the concentration inequalities of the mean-field interaction particle system (\ref{MFS}), as applications of   our main theorems.
This kind of concentration estimate are useful to numerical simulations and mean-field limit. Under the conditions that $V$ is uniformly
convex and $W$ is convex, Malrieu \cite{[MF1]} established logarithmic
Sobolev inequality and then used its connection with optimal transport and concentration of measure to get the following non-asymptotic
bounds on the deviation of the empirical mean of an observable $f$ from its true mean,
\begin{equation}\label{Mal}
\sup_{\|f\|_{Lip}\le 1}\mathbb{P}\left\{|\frac1N\sum_{i=1}^N f(X_t^{i,N}) - \mu_t(f)| > \frac{A}{\sqrt{N}} + \delta
 \right\} \le 2 e^{-\lambda N \delta^2},\ t>0,\ \delta\ge 0
\end{equation}
where $A$ and $\lambda$ are positive constants depending on the particle system.

As pointed out in \cite{[BGV]}, this approach can lead to nice bounds but it is limited to a finite number of observables. Bolley-Guillin-Villani
\cite[Theorem 2.9]{[BGV]} obtained for any $t>0$ fixed and $\delta>0$
\begin{equation}\label{BGVc}
\mathbb{P}\left\{\sup_{\|f\|_{Lip}\le 1}|\frac1N\sum_{i=1}^N f(X_t^{i,N}) - \mu_t(f)| > \delta
 \right\} \le C(1+t\delta^{-2}) e^{-K(t) N \delta^2},
\end{equation}
for all $N$ big enough (quantifiable), where $K(t)$ depending on $t$ is some explicitly computable constant.
Furthermore, Bolley \cite{[B]} got quantitative concentration inequalities on the sample path space with uniform norm, on a given time interval $[0,T]$,
which implies (\ref{BGVc}) by projection at time $t\in [0,T]$.

\subsection{Uniform in time concentration inequality} Our previous general results allow us to generalize (\ref{Mal}) and to reinforce (\ref{BGVc}) (under stronger conditions of course).

\bprop\label{prop51} Assume {\bf (H)}, (\ref{b1}) and (\ref{b5}). Let $X_t^{(N)}=(X_t^{1,N},\cdots,X_t^{N,N}),\forall t\ge 0,$
then for any Lipschitzian function $F$ on $(\rr^d)^N$, we have for any lower bounded convex function $\varphi$ on $\rr$,
\bbeq\label{prop51a}
\ee_x \varphi\left(F(X^{(N)}_T)- \ee_x F(X^{(N)}_T)\right) \le \ee \varphi\left( \alpha A_\vep \sqrt{\frac{N}{2K_\vep}}   \xi \right),\ \forall x\in(\rr^d)^N,\ \forall T>0
\nneq
where $\xi$ is some standard real Gaussian random variable of law $\NN(0,1)$, $\alpha:=\|F\|_{Lip(d_{l^1})}=\max_{1\le i \le N}\|\nabla_iF\|_\infty$,   $A_{\vep}$ and $K_{\vep}$ are  given in Theorem \ref{W1c}.

In particular for any initial distribution $\mu_0$ satisfying the Gaussian integrability assumption on $\rr^d$, we have for any $\delta, T>0$
\begin{equation}\label{prop51b}
\mathbb{P}_{\mu_0^{\otimes N}}\left\{F(X^{(N)}_T)-\mathbb{E}_{\mu_0^{\otimes N}} F(X_{T}^{(N)})>\delta \right\} \le \exp\left(- \frac{K_\vep\delta^2}{N \alpha^2 A_{\vep}^2 \left[1+ 2 c_G(\mu_0) K_\vep e^{-2K_\vep T}\right]}\right).
\end{equation}

\nprop

\bprf Without loss of generality we may  assume that $\alpha=\max_{1\le i \le N}\|\nabla_iF\|_\infty=1$.

By approximation we may assume that $F$ is $C^2$-smooth with bounded derivatives of the first and the second order.
For any initial position $x\in(\rr^d)^N$, let $M_t=\mathbb{E}_x(F(X^{(N)}_{T})|\mathcal{F}_t),\  0 \le t \le T$.
Then by applying It\^o's formula to $u(t,x)=P_{T-t}F(x)$, we have
\begin{equation}\label{EX2a}
F(X^{(N)}_{T})- \mathbb{E}_x F(X^{(N)}_{T})= M_T-M_0=\sum_{i=1}^N \int_0^T \nabla_{i} P_{T-t}F(X^{(N)}_{t}) dB_t^i,
\end{equation}

Note that by Theorem \ref{W1c}, for any $\vep>0$ such that $K_\vep>0$,
we have
\begin{equation}\label{EX2b}
W_{d_{l^1}}(P_t^{(N)}(x,\cdot), P_t^{(N)}(y,\cdot)) \le A_{\vep} e^{-K_{\vep}t} d_{l^1}(x, y), \forall x,y \in (\mathbb{R}^d)^N,
\end{equation}
which implies that
\begin{equation}\label{EX2c}
|\nabla_{i} P_{T-t}F|\le A_{\vep} e^{-K_{\vep}(T-t)}, \ 1\le i \le N,
\end{equation}
where $A_{\vep}$ and $K_{\vep}$ are the same as given in Theorem \ref{W1c}.

Since $M_t=\xi_{\tau_t}$ where $(\xi_t)$ is a real valued Brownian motion w.r.t. some new filtration $(\tilde \FF_t)$ and
$$\tau_t=\<M\>_t=\int_0^t\sum_{i=1}^N |\nabla_{i} P_{t-s}F(X^{(N)}_{s})|^2 ds\le
 \frac{A_{\vep}^2}{2K_{\vep}} N=:CN$$ is a stopping time w.r.t. $(\tilde \FF_t)$,  we obtain
\begin{equation}\label{C4}
\begin{aligned}
\mathbb{E}_x \varphi\left(F(X^{(N)}_{T})-\mathbb{E}_x F(X^{(N)}_{T})\right) &= \mathbb{E} \varphi\left(M_T-M_0\right) = \ee \varphi(\xi_{\tau_T}) \\
&=\ee\varphi\left(\ee (\xi_{CN}| \tilde \FF_{\tau_T})\right) \\
&\le \ee \varphi\left( \xi_{CN}\right)\ \ \ \ \ \text{(by Jensen's inequality) }\\
&=\ee \varphi\left( A_\vep \sqrt{\frac{N}{2K_\vep}}   \xi_1 \right)\\
\end{aligned}\end{equation}
the desired result (\ref{prop51a}).

Letting $g(x):=\ee_{x}F(X_T^{(N)}),\forall x\in(\rr^d)^N$. By (\ref{EX2c}) we have \bbeq\label{C6}\|g\|_{Lip(d_{l^1})}=\max_{1\le i\le N}\|\nabla_i g\|_\infty\le A_\vep e^{-K_\vep T}.\nneq
Applying (\ref{prop51a}) to $\varphi(z)=e^{\lambda z}$ ($\lambda\in \rr$), we get
$$
\aligned
&\ee_{\mu_0^{\otimes N}} \exp\left(\lambda [F(X_T^{(N)})-\ee_{\mu_0^{\otimes N}}F(X_T^{(N)})] )\right)\\
&=\int_{(\rr^d)^N} \ee_x \exp\left(\lambda [F(X_T^{(N)})-\ee_{x}F(X_T^{(N)})]\right) \cdot \exp\left(\lambda [g(x)-\mu_0^{\otimes N}(g)]\right)d\mu_0^{\otimes N}(x)\\
&\le \int_{(\rr^d)^N} \ee \exp\left(\lambda A_\vep \sqrt{\frac N{2K_\vep}}\xi_1\right) \cdot \exp\left(\lambda [g(x)-\mu_0^{\otimes N}(g)]\right)d\mu_0^{\otimes N}(x)\\
&\le \exp\left(\frac{NA_\vep^2\lambda^2}{4K_\vep}\right)\exp\left(\frac{\lambda^2}{2} N c_G(\mu_0)\|g\|^2_{Lip(d_{l^1})} \right)\\
&\le  \exp\left(\frac{N\lambda^2 A_\vep^2}{2}\left[\frac 1{2K_\vep}+ c_G(\mu_0) e^{-2K_\vep T}\right]\right)
\endaligned
$$
where the third and the last inequality follows from the Gaussian concentration condition on the initial distribution $\mu_0$ (see \eqref{Uc} in Remark \ref{remc}) and \eqref{C6} respectively.

Finally the concentration inequality (\ref{prop51b}) is derived from the above inequality by the standard procedure of Chebyshev's inequality and optimization over $\lambda$.
\nprf

\begin{example}\label{EX2} {\rm Given a Lipschitzian observable $f:\rr^d\to\rr$ with $\|f\|_{Lip}=1$ and $N\ge2$, let $F(x)=\frac 1N \sum_{i=1}^N f(x^i)$.
For any $T>0$,
$$F(X^{(N)}_{T})=\frac{1}{N}\sum_{i=1}^N  f(X^{i,N}_T)$$
is the empirical mean of $f$ at time $T$.
Since $$\|F\|_{Lip_{(d_{l^1})}} = \frac{1}{N} \|f\|_{Lip}=\frac 1N$$
 we obtain by (\ref{prop51b}) for any $\delta>0$,

\begin{equation}\label{prop51b}
\mathbb{P}_{\mu_0^{\otimes N}}\left\{\frac{1}{N}\sum_{i=1}^N  f(X^{i,N}_T)-\mathbb{E}_{\mu_0^{\otimes N}} f(X_{T}^{1,N})>\delta \right\} \le \exp\left(- \frac{N K_\vep\delta^2}{  A_{\vep}^2 \left[1+ 2 c_G(\mu_0) K_\vep e^{-2K_\vep T}\right]}\right).
\end{equation}
As the absolute value of the bias $|\mathbb{E}_{\mu_0^{\otimes N}} f(X_{T}^{1,N}) - \mu_T(f)| \le A/\sqrt{N}$ by Remark \ref{Prop43Rem2}, our result above generalizes Malrieu's result (\ref{Mal}) to the case that $V$ may have many wells.

}
\end{example}

\subsection{Concentration for time average}
The counterpart of Proposition \ref{prop51} for the empirical time average is presented in the following result.

\bprop\label{C} Assume  {\bf (H)}, (\ref{b1}) and (\ref{b2}). Given any $T\in(0,+\infty]$,
let $F$ be any $d_{L^1[0,T]}$-Lipschitzian continuous function on $C([0,T], (\mathbb{R}^d)^N)$, given by
$$F(X^{(N)}_{[0,T]}):=G\left(\int_0^T g_1(X^{(N)}_t)dt, \cdots, \int_0^T g_n(X^{(N)}_t)dt\right),$$
where $G\in C^2(\rr^n)$, $g_i \in C^2((\mathbb{R}^d)^N, \mathbb{R}),\ 1\le i \le n$. Then for any convex function $\varphi$ on $\mathbb{R}$ and any starting point $X_0^{(N)}=x\in(\rr^d)^N$,
we have
\begin{equation}\label{C1}
\mathbb{E}_x \varphi\left(F(X^{(N)}_{[0,T]})-\mathbb{E}_x F(X^{(N)}_{[0,T]})\right) \le \ee \varphi\left(\sqrt{ NT}
\|F\|_{Lip(d_{L^1[0,T]})} c_{Lip}\xi\right),
\end{equation}
where  $\xi$ is some standard real Gaussian random variable of law $\NN(0,1)$, and
$$c_{Lip}=\frac{h^{\prime}(0)}{1-\|\nabla^2_{xy}W\|_{\infty} \|h^{\prime}\|_{\infty}}.$$
\nprop

\bprf Let $\{\mathcal{F}_t\}_{t\ge 0}$ be the filtration generated by the process $(X^{(N)}_t)_{t\ge 0}$
and $$M_t=\mathbb{E}(F(X^{(N)}_{[0,T]})|\mathcal{F}_t),\  0 \le t \le T .$$ Then by the martingale representation theorem, we have
\begin{equation}\label{C2}
F(X^{(N)}_{[0,T]})- \mathbb{E} F(X^{(N)}_{[0,T]})= M_T-M_0=\sum_{i=1}^N \int_0^T \beta_t^i dB_t^i,
\end{equation}
where $\beta_t^i, 1\le i \le N$ are adapted processes w.r.t. $\mathcal{F}_t$, and $B_t^i, 1\le i \le N$ are $N$ independent standard Brownian motions on $\rr^d$.

Let $A_t^k=\int_0^t g_k(X^{(N)}_s)ds,\ 1\le k\le n,$ and $A_t=(A_t^1, \cdots, A_t^n)$. Note that
$$M_t=\phi(A_t, X^{(N)}_t)$$
where
$$
\phi(a, x):=\ee \left( G\left(a_1 +\int_t^T g_1(X^{(N)}_s)ds, \cdots, a_n +\int_t^T g_N(X^{(N)}_s)ds\right)|X^{(N)}_t =x\right), $$
for $a=(a_1,\cdots, a_n)\in \rr^n, x\in (\rr^d)^N$.  Since $\varphi$ is $C^2$ (for $V, W$ are $C^2$), we can
apply It\^o's formula to obtain that
$$\beta_t^i=\partial_{x_i} \varphi (A_t, X^{(N)}_{t}).$$

For any $x=(x^1, \cdots, x^i, \cdots, x^N)\in (\mathbb{R}^{d})^N$, denote $y=(x^1, \cdots, y^i, \cdots, x^N)$
which only differs from $x$ at the $i$-th coordinate. Let $(X_t^{(N)})_{t\ge0}, (Y^{(N)}_t)_{t\ge0}$ be an optimal coupling of $\pp_x, \pp_y$ for
$W_{1,d_{L^1[0,T]}}(\pp_x,\pp_y)$ (this optimal coupling exists because $d_{L^1[0,T]}$ is lower semi-continuous from $(C(\rr^+,(\rr^d)^N))^2$ to $[0,+\infty]$).
Then for any $0\le t \le T$ and $i=1, \cdots, N,$ we have
\begin{equation}\label{C3}
\begin{aligned}
&|\partial_{x_i} \phi(a,x)| \le \limsup_{y^i \rightarrow x^i}\frac{|\phi(a,x)-\phi(a,y)|}{|x^i-y^i|}|\\
 &=\limsup_{y^i \rightarrow x^i} \frac{1}{|x^i-y^i|}|\ee[ G(a+\int_0^{T-t} g(X_s^{(N)}ds)] - \ee[ G(a+\int_0^{T-t} g(Y_s^{(N)}ds)]|\\
 & \le  \limsup_{y^i \rightarrow x^i}  \frac{\|F\|_{Lip(d_{L^1[0,T]})}}{|x^i-y^i|}\ee \int_0^{\infty} d_{l^1}(X^{(N)}_s, Y^{(N)}_s) ds\\
 &= \|F\|_{Lip(d_{L^1[0,T]})} \limsup_{y^i \rightarrow x^i}  \frac{W_{1, d_{L^1}}(\pp_x, \pp_y)}{|x^i-y^i|}\\
 &\le \|F\|_{Lip(d_{L^1[0,T]})}\cdot c_{Lip}
\end{aligned}\end{equation}
where the last inequality follows by Theorem \ref{W1}.

We now repeat the argument in the proof of Proposition \ref{prop51}.
Since $\sum_{i=1}^N\int_0^T \beta_t^i dB_t^i=\xi_{\tau_T}$ where $(\xi_t)$ is a real valued Brownian motion w.r.t. some new filtration $(\tilde \FF_t)$ and $\tau_T=\int_0^T\sum_{i=1}^N |\beta_t^i|^2 dt\le \|F\|_{Lip(d_{L^1[0,T]})}^2 c_{Lip}^2 NT=:CNT$ is a stopping time w.r.t. $(\tilde \FF_t)$, we obtain
\begin{equation}\label{C7}
\begin{aligned}
\mathbb{E}_x \varphi\left(F(X^{(N)}_{[0,T]})-\mathbb{E} F(X^{(N)}_{[0,T]})\right) &= \mathbb{E} \varphi\left(\sum_{i=1}^N
\int_0^T \beta_t^i dB_t^i\right) = \ee \varphi(\xi_{\tau_T}) \\
&=\ee\varphi\left(\ee (\xi_{CNT}| \tilde \FF_{\tau_T})\right) \\
&\le \ee \varphi\left( \xi_{CNT}\right)\ \ \ \ \ \text{(by Jensen's inequality) }\\
&=\ee \varphi\left( \sqrt{NT} \|F\|_{Lip(d_{L^1[0,T]})} c_{Lip} \xi_1 \right)\\
\end{aligned}\end{equation}
the desired result.
\nprf

Next we give the proof of Corollary \ref{U}.
\bprf[Proof of Corollary \ref{U}] For any given $\lambda, T>0$, let
$$F(X_{[0,T]}^{(N)})=\frac{1}{T}\int_0^T U_N(f_m)(X_t^{(N)})dt. $$
Since $f_m$ is $1$-Lipschitzian w.r.t. the $d_{l^1}$-metric on $(\rr^d)^m$, by an easy calculation we have
$$\|F\|_{Lip_{(d_{L^1[0,T]})}} \le \frac{m}{NT}.$$

Let $g(x)=\ee_x F,\ \forall x\in(\rr^d)^N$. For any fixed initial value $x\in(\rr^d)^N$, by applying Proposition \ref{C} with $\varphi(z)=e^{\lambda z}$,  we get
\begin{equation}\begin{aligned}\label{cor3a}
&\mathbb{E}_x\exp\left(\lambda \left[\frac 1T\int_0^T U_N(f_m)(X_t^{(N)}) dt - g(x)\right]\right)\\
\le & \ee \exp\left(\frac{m\lambda}{\sqrt{NT}}c_{Lip}\xi_1\right)=  \exp \left(\frac{m^2\lambda^2c_{Lip}^2}{2NT}\right).\\
\end{aligned}\end{equation}
 By the proof of Proposition \ref{C},
$$\|g\|_{Lip(d_{l^1})}\le c_{Lip} \|F\|_{d_{L^1[0,T]}}\le \frac{m c_{Lip}}{NT}.$$
By the condition (\ref{Ub}) and its consequence (\ref{Uc}), the product measure $\mu_0^{\otimes N}$ satisfies
\bbeq\label{cor3e}\aligned
\int_{(\rr^d)^N} e^{\lambda (g - \mu_0^{\otimes N}(g))} d\mu_0^{\otimes N} &\le \exp\left(\frac{1}{2}N c_G(\mu_0) \lambda^2\|g\|^2_{Lip(d_{l^1})}\right)\\
&\le\exp\left(\frac{1}{2NT^2}c_G(\mu_0) \lambda^2 m^2 c_{Lip}^2\right) .
\endaligned\nneq
Hence for the i.i.d. initial values $X_0^{1,N},\cdots, X_0^{N,N}$ with the common law $\mu_0$, noting that
$$
\ee \frac 1T \int_0^T U_N(f_m)(X_t^{(N)}) dt = \mu_0^{\otimes N}(g)
$$
we have
 \begin{equation}\begin{aligned}\label{cor3b}
 &\ee \exp\left(\lambda \left[\frac 1T \int_0^T U_N(f_m)(X_t^{(N)}) dt -\ee \frac 1T \int_0^T U_N(f_m)(X_t^{(N)}) dt\right]\right)\\
 =& \int_{(\rr^d)^N}\ee_x\left[\exp\left(\lambda \left[ \frac 1T \int_0^T U_N(f_m)(X_t^{(N)}) dt - g(x)\right]\right)\right] e^{\lambda(g(x)-\mu_0^{\otimes N}(g) )}d\mu_0^{\otimes N}(x)\\
 \le & \exp\left(\frac{m^2\lambda^2c_{Lip}^2}{2NT}\right) \int_{(\rr^d)^N} e^{\lambda(g(x)-\mu_0^{\otimes N}(g) )}d\mu_0^{\otimes N}(x) \\
 \le & \exp\left(\frac{m^2\lambda^2c_{Lip}^2}{2NT}\left(1+\frac{c_G(\mu_0)}{T}\right)\right)
 \end{aligned}\end{equation}
 where the second  inequality follows from (\ref{cor3a}), and the last inequality is a consequence of (\ref{cor3e}). This gives us (\ref{cor3}). Finally (\ref{cor4}) follows from  (\ref{cor3}), by the standard procedure of Chebyshev's inequality and optimization over $\lambda>0$.
\nprf

\brmk{\rm The time-particle average $\frac 1{NT}\int_0^T f(X^{i,N}_t)dt$ is used to approximate $\mu_\infty(f)$ where $\mu_\infty=\lim_{t\to+\infty}\mu_t$ is the unique equilibrium state of the McKean-Vlasov equation (a consequence of Corollary \ref{cor11} by Banach's fixed point theorem). For applying Corollary \ref{U}, it remains to bound the bias
$$\aligned
&|\ee_{\mu_0^{\otimes N}} \frac 1{NT}\int_0^T \sum_{i=1}^N f(X^{i,N}_t)dt-\mu_\infty(f)|\\
&\le |\ee_{\mu_0^{\otimes N}} \frac 1{T}\int_0^T [f(X^{1,N}_t)-\mu_t(f)]dt| + \frac 1T \int_0^T|\mu_t(f)-\mu_\infty(f)| dt \\
&\le  \frac 1{T}\int_0^T |\mu_t^{1,N}(f)-\mu_t(f)|dt+ \|f\|_{Lip}\frac 1T \int_0^TW_1(\mu_t,\mu_\infty) dt \\
&\le  \|f\|_{Lip} \left(\sup_{t\ge 0}W_1(\mu_t^{1,N},\mu_t) + \frac 1T \int_0^TW_1(\mu_t,\mu_\infty) dt\right) \\
&\le \|f\|_{Lip} \left(\frac {A}{\sqrt{N}} + \frac {B}{T}\right)
\endaligned$$
where in the last inequality, the first term comes from the uniform in time propagation of chaos (\ref{CP0}) in Theorem \ref{Wp}, and the second follows by (\ref{cor11a}) in Corollary \ref{cor11}. We believe that the bias should be of order $O(1/N +1/T)$, but we do not know how to prove it.
}\nrmk

\bibliographystyle{plain}

\end{document}